\documentclass[11pt]{article}

\usepackage{amscd,amstex,amssymb}

\def\eqref#1{(\ref{#1})}
\newcommand{\goth}{\mathfrak}
\newcommand{\arrow}{{\:\longrightarrow\:}}
\newcommand{\Z}{{\Bbb Z}}
\newcommand{\C}{{\Bbb C}}
\newcommand{\R}{{\Bbb R}}

\newcommand{\1}{\sqrt{-1}\:}
\newcommand{\inangles}[1]{{\langle #1\rangle}}

\newcommand{\restrict}[1]{{\left|_{{\phantom{|}\!\!}_{#1}}\right.}}

\renewcommand{\c}[1]{{\cal #1}}
\newcommand{\calo}{{\cal O}}

\renewcommand{\tilde}{\widetilde}
\renewcommand{\bar}{\overline}
\renewcommand{\phi}{\varphi}
\renewcommand{\epsilon}{\varepsilon}

\newcommand{\im}{\operatorname{im}}
\newcommand{\End}{\operatorname{End}}

\newcommand{\Sec}{\operatorname{Sec}}
\newcommand{\Hor}{\operatorname{Hor}}
\newcommand{\Def}{\operatorname{Def}}
\newcommand{\Tw}{\operatorname{Tw}}

\newcommand{\comment}[1]{{}}

\def\blacksquare{\hbox{\vrule width 4pt height 4pt depth 0pt}}
\def\endproof{\blacksquare}

\makeatletter

\@ifundefined{Bbb}
     {\newcommand{\Bbb}[1]{{\mathbb #1}}}%
{}

\newcommand{\ps@verbit}{%
  \renewcommand{\@oddhead}{%
          \scriptsize
          {Misha Verbitsky \ \ Hypercomplex varieties}
          \hfil\tiny {final version, February 17, 1997}}
  \renewcommand{\@evenhead}{\@oddhead}
  \renewcommand{\@oddfoot}{\hfil\thepage\hfil}
  \renewcommand{\@evenfoot}{\@oddfoot}}
 
\newcounter{Mycounter}[section]
\newcounter{lemma}[section]
\renewcommand{\thelemma}{{Lemma \thesection.\arabic{lemma}}}
\newcommand{\lemma}{%
     \setcounter{lemma}{\value{Mycounter}}
     \refstepcounter{lemma}
     \stepcounter{Mycounter}
     {\bf \thelemma:\ }}

\newcounter{claim}[section]
\renewcommand{\theclaim}{{Claim \thesection.\arabic{claim}}}
\newcommand{\claim}{%
     \setcounter{claim}{\value{Mycounter}}
     \refstepcounter{claim}
     \stepcounter{Mycounter}
     {\bf \theclaim:\ }}

\newcounter{sublemma}[section]
\newcounter{corollary}[section]
\newcounter{theorem}[section]
\renewcommand{\thetheorem}{{Theorem \thesection.\arabic{theorem}}}
\newcommand{\theorem}{%
     \setcounter{theorem}{\value{Mycounter}}
     \refstepcounter{theorem}
     \stepcounter{Mycounter}
     {\bf \thetheorem:\ }}

\newcounter{conjecture}[section]
\newcounter{proposition}[section]
\renewcommand{\theproposition}
       {{Proposition \thesection.\arabic{proposition}}}
\newcommand{\proposition}{%
     \setcounter{proposition}{\value{Mycounter}}
     \refstepcounter{proposition}
     \stepcounter{Mycounter}
     {\bf \theproposition:\ }}

\newcounter{definition}[section]
\renewcommand{\thedefinition}
       {{Definition \thesection.\arabic{definition}}}
\newcommand{\definition}{%
     \setcounter{definition}{\value{Mycounter}}
     \refstepcounter{definition}
     \stepcounter{Mycounter}
     {\bf \thedefinition:\ }}

\newcounter{example}[section]
\newcounter{remark}[section]
\renewcommand{\theremark}{{Remark \thesection.\arabic{remark}}}
\newcommand{\remark}{%
     \setcounter{remark}{\value{Mycounter}}
     \refstepcounter{remark}
     \stepcounter{Mycounter}
     {\bf \theremark:\ }}

\newcounter{problem}[section]
\newcounter{question}[section]
\@addtoreset{equation}{section}
\@addtoreset{footnote}{section}
\makeatother

\begin{document}

\begin{center}
{\Large\bf
Hypercomplex varieties}\\[4mm]
Misha Verbitsky,\footnote{Supported by the NSF grant 9304580}\\[4mm]
{\tt verbit@@thelema.dnttm.rssi.ru, verbit@@math.ias.edu}
\end{center}

\hfill

{\small 
\hspace{0.2\linewidth}
\begin{minipage}[t]{0.7\linewidth}
We give a number of equivalent definitions of
hypercomplex varieties and construct a twistor
space for a hypercomplex variety. We prove that 
our definition of a hypercomplex variety
(used, e. g., in alg-geom 9612013) is equivalent to a definition 
proposed by Deligne and Simpson, who used twistor spaces.
This gives a way to define hypercomplex spaces
(to allow nilpotents in the structure sheaf).
We give a self-contained proof of desingularization theorem for
hypercomplex varieties: a normalization of a hypercomplex
variety is smooth and hypercomplex. 
\end{minipage}
}

\tableofcontents


\section{Introduction}


Hypercomplex varieties are a very natural class of objects. 
This notion allows one to speak uniformly of a number of 
disparate examples coming from hyperk\"ahler geometry and the
theory of moduli spaces (\ref{_triana_hyperco_Remark_}, Subsection 
\ref{_hyperholomo_Subsection_}). 

A main property of hypercomplex varieties is that they can
be desingularized in a functorial way, and this desingularization
is a hypercomplex manifold. Moreover, this desingularization is
achieved by taking normalization. This is a very useful result --
see \cite{_Arapura_}, \cite{_Verbitsky:New_} for some of its uses.

This is why it is important to dwell on the definition of hypercomplex
varieties. We give a number of equivalent definitions, in order to produce 
conditions which are easy to check in every situation.

It is now possible to prove that something is {\bf not} hypercomplex.
In \ref{_quotie_not_hyperco_Proposition_}, 
we show that for a hypercomplex variety $M$, and a finite
group $G$ acting on $M$ by automorphisms, the quotient
$M/G$ is {\bf not} hypercomplex, unless $G$ acts on $M$
freely. More precisely, we show that, for $M_0\subset M$ 
the part of $M$ where $G$ acts freely, the natural hypercomplex
structure on $M_0/G$ cannot be extended to $M/G$. 

\subsection{An overview}

The definitions of a hypercomplex variety,
given in  \cite{_Verbitsky:Desingu_} and 
\cite{_Verbitsky:DesinguII_}, on one hand, and in
\cite{_Verbitsky:Deforma_}, \cite{_Verbitsky:Hyperholo_bundles_} 
on the other hand, 
are not identical. We show that these two definitions
are in fact equivalent.

We give the weakest form of the definition
of hypercomplex varieties (one used in
\cite{_Verbitsky:Deforma_}, \cite{_Verbitsky:Hyperholo_bundles_}), 
and give the proof of the desingularization theorem under these 
assumptions. This proof is essentially identical to the proof used
in \cite{_Verbitsky:Desingu_} and \cite{_Verbitsky:DesinguII_},
though our assumptions are weaker, and the
arguments are more rigorous. Then we obtain the
stronger version of the definition 
(used in \cite{_Verbitsky:Desingu_} and 
\cite{_Verbitsky:DesinguII_}) from the 
desingularization and a recent result of Kaledin
\cite{_Kaledin_}.

An almost hypercomplex variety (\ref{_almost_hyperco_Definition_})
is a real analytic variety $M$ with a quaternionic action on its
sheaf of real analytic differentials.
For each $L\in \Bbb H$, $L^2 = -1$, $L$ gives 
an almost complex structure on $M$. This almost
complex structure is called {\bf induced by the quaternionic
action}. We say that $M$ is {\bf hypercomplex} 
(\ref{_hyperco_Definition_}) if there are 
induced almost complex structures $I, J \in \Bbb H$, such that
$I\neq \pm J$, and $I$, $J$ are integrable 
(\ref{_integra_almost_comple_Definition_}).
We show that then {\bf all} induced complex structures
are integrable. This was {\em a definition} of the 
hypercomplex variety which we used in \cite{_Verbitsky:Desingu_} and 
\cite{_Verbitsky:DesinguII_}. 

For almost hypercomplex manifolds
(i. e., smooth almost hypercomplex varieties),
D. Kaledin \cite{_Kaledin_} proved that integrability
of two $I, J \in \Bbb H$, $I\neq \pm J$ implies integrability
of all induced complex structures. We prove 
the desingularization theorem (\ref{_desingu_Theorem_}) 
for hypercomplex varieties, and then apply Kaledin's result
to obtain integrability of all induced complex structures
(\ref{_all_indu_comple_integra_Theorem_}).

The last part of this paper deals with several other versions
of the definition of hypercomplex variety, which are,
as we show, equivalent.

We define the {\bf twistor space} of a hypercomplex
variety (\ref{_twistor_Definition_}), constructed
as an almost
complex variety. Using Kaledin's theorem and
desingularization, we prove that the almost complex
structure on the twistors is integrable, i. e., the twistor
space is a complex variety. We give a description of the
hypercomplex structure in terms of the twistor space, 
 following \cite{_HKLR_}, \cite{_Simpson:hyperka-defi_} 
and \cite{_Deligne:defi_}. 

The twistor space $\Tw$ is
a complex variety equipped with a holomorphic
map $\pi:\; \Tw \arrow \C P^1$ and an anticomplex
involution $\iota:\; \Tw \arrow \Tw$. The original hypercomplex
variety is identified with a fiber $\pi^{-1}(I)$, $I\in \C P^1$,
and the data $(\Tw, \pi, \iota)$ are essentially sufficient
to recover the hypercomplex structure on $M = \pi^{-1}(I)$
(see \eqref{_twistor_data_Equation_} 
for details and a precise statement).

It is possible to define a hypercomplex variety in terms
of $(\Tw, \pi, \iota)$. This definition was proposed
by Deligne and Simpson (\cite{_Deligne:defi_}, 
\cite{_Simpson:hyperka-defi_}). We show that their
definition is equivalent to ours. 
Their definition has the significant advantage
that it does not assume that $M$ is reduced,
and indeed might be used to define hypercomplex
spaces, i. e., to allow nilpotents in the 
structure sheaf.

\subsection{Twistor spaces: an introduction}

The twistor space of a hypercomplex variety is the following
object. 

Let $M$ be a hypercomplex variety. Then the quaternion algebra
$\Bbb H$ acts in $\Omega^1M$ in such a way that for all 
$L\in \Bbb H$, $L^2 =-1$, the corresponding operator
is an almost complex structure. Identifying the set
$\{ L\in \Bbb H \;\; |\;\; L^2=-1\}$ with $\C P^1$, we obtain 
an almost complex structure $\c I$ on $\C P^1 \times M$ acting on
$T _{L,m} \C P^1 \times M = T_L \C P^1 \oplus T_m M$
as $I_{\C P^1} \oplus L$, where $L\in \C P^1$ acts on $T_m M$
as the corresponding quaternion, and $I_{\C P^1}$ is the complex
structure operator on $\C P^1$. This almost complex structure
is integrable (\ref{_twi_integra_Claim_}; essentially, this is 
proven by D. Kaledin \cite{_Kaledin_}). The corresponding
variety is called {\bf the twistor space of $M$}, denoted by $\Tw$ 
(\ref{_twistor_Definition_}).

Consider the holomorphic projection map $\pi:\; \Tw\arrow \C P^1$.
Let $\iota_0:\; \C P^1 \arrow \C P^1$ be the anticomplex involution 
with no fixed points given by the central symmetry of $S^2 \subset \R^3$,
and $\iota:\; \Tw \arrow \Tw$ map $(s, m)$ to $(\iota_0(s), m)$.
Clearly, $\iota$ is also an anticomplex involution. 

The set of sections of the projection $\pi$ is called
{\bf the space of twistor lines}, denoted by $\Sec$.
This space is equipped with complex structure,
by Douady (\cite{_Douady_}). 
 
Consider a twistor line 
$I \stackrel {s_m} \arrow (I \times m)\in \C P^1 \times M = \Tw$,
where $m\in M$. Then $s_m$ is called {\bf a horisontal twistor line}.
The variety of horisontal twistor lines is denoted by 
$\Hor$. Clearly, the line $\im(s_m)\subset \Tw$ is fixed 
by the involution $\iota$. One can show that a space of 
horisontal twistor lines is a
connected component of the space $\Sec^\iota$ of all 
twistor lines fixed by $\iota$. Clearly, the real analytic
variety $\Hor$ is identified with $M$, and the induced complex
structures on $M$ come from the natural 
isomorphisms between $\Hor$ and fibers of $\pi$.
This shows how to recover
the hypercomplex variety $M$ from $\Tw$, $\iota$ and $\pi$.

\hfill

According to Deligne and Simpson (\cite{_Deligne:defi_}, 
\cite{_Simpson:hyperka-defi_}), singular hyperk\"ahler varieties
should be defined in terms of the following 
data (see \eqref{_twistor_data_Equation_} for details).

\begin{description}
\item[(i)] $\Tw, \pi, \iota$
\item[(ii)] a component $\Hor$ of the variety of
$\iota$-invariant holomorphic sections of $\pi$.
\end{description}
satisfying the following conditions 
(see \eqref{_twistor_properties_Equation_} for details).
\begin{description}
\item[(a)] Through every point $x\in \Tw$ passes a unique
line $s\in \Hor$
\item[(b)] Let $s\in Hor$. Consider the points $x$, $y$ situated
in a small neighbourhood $U$ of $\im s \subset \Tw$, $\pi(x)\neq \pi(y)$.
Assume that $x$ and $y$ belong to the same irreducible component of $U$.
Then there exist a unique twistor line $s_{xy}\in \Sec$ passing
through $x$, $y$.
\end{description}
The advantage of the definition of Deligne--Simpson
(see \ref{_Deligne_Simpson_Definition_})
is that it is easy to adopt for the scheme situation: 
one may allow nilpotents in
structure sheaf. For the definition of a {\bf hypercomplex space},
see \ref{_hyperco_spaces_Definition_}.

We show that Deligne--Simpson's definition of a hypercomplex variety
is equivalent to ours (\ref{_hype_vari_and_twi_equiva_Theorem_}).
The condition (b) in the above listing is very difficult to check.
We replace it by
an equivalent condition, which should be thought of as its
linearization (see \ref{_twi_hyperco_type_Definition_} for details).

\begin{description}
\item[(b$'$)] For each $s\in \Hor$, the conormal sheaf
\[ \ker \bigg (\Omega^1(\Tw)\restrict{\im s} \arrow \Omega^1(\im s)\bigg) \]
is isomorphic to $\oplus(\calo(-1))$.
\end{description}

\ref{_Deli_Simpsi_equi_infinite_Theorem_} shows that
the above conditions (b) and (b$'$) are equivalent.

The condition (b$'$) is very easy to check, and hence is
extremely useful. In Subsection
\ref{_hyperholomo_Subsection_}, we use the 
condition (b) to give a new proof
that a deformation  space of a stable holomorphic 
bundle $B$ over a hyperk\"ahler manifold is 
a hyperk\"ahler variety, provided
that the Chern classes of $B$ are $SU(2)$-invariant.

\subsection{Desingularization of hypercomplex varieties}

Let $M$ be a hypercomplex variety (\ref{_hyperco_Definition_}),
$I$ an integrable induced complex structure. Consider 
the complex variety $(M,I)$, which is $M$ with  the
complex structure defined by $I$. The Desingularization Theorem
(\ref{_desingu_Theorem_}) says that the normalization
$\widetilde{(M,I)}$ of $(M,I)$ is smooth and hypercomplex.  
This desingularization is canonical and functorial: the hypercomplex
manifold $\widetilde{(M,I)}$ is independent from the choice of 
induced complex structure $I$.

The proof of the Desingularization Theorem goes along 
the following lines. We define spaces with locally 
homogeneous singularities (\ref{_SLHS_Definition_}).
A  space with locally homogeneous singularities 
(SLHS) is an analytic space $X$ such that for all $x\in X$, the
$x$-completion of a local ring $\calo_xX$ is isomorphic
to an $x$-completion of the associated graded ring 
$(\calo_xX)_{gr}$. We show that hypercomplex varieties are always
SLHS (\ref{_hyperco_SLHS_Theorem_}). This is proven using 
an elementary argument from commutative algebra.

We work with a complete local Noetherian
ring $A$ over $\C$, with a residual field $\C$. 
By definition, an automorphism $e:\; A \arrow A$
is called {\bf homogenizing} (\ref{_homogeni_automo_Definition_})
if its differential acts as
a dilatation on the Zariski tangent space of $A$, with dilatation
coefficient $|\lambda|<1$. As usual, by the Zariski tangent space
we understand the space $(\goth m_A /\goth m_A^2)^*$, where $\goth m_A$
is a maximal ideal of $A$. For a ring $A$
equipped with a homogenizing automorphism
$e:\; A \arrow A$, we show that $A$ has locally homogeneous 
singularities.

We complete the proof that all hypercomplex varieties are 
SLHS by constructing an explicit homogenizing automorphism in a
local ring of germs of holomorphic functions on a complex
variety underlying a given hypercomplex variety 
(\ref{_homogenizing_Proposition_}).

The proof of Desingularization Theorem proceeds then with
a study of a tangent cone of a hypercomplex variety. 
For every point of a hypercomplex variety, the corresponding
Zariski tangent space $T_xM$ is equipped with a quaternionic action.
This makes $T_zM$ into a hyperk\"ahler manifold, with appropriately
chosen metric. We show that the tangent cone $Z_xM$ of $M$, 
considered as a subvariety of $T_xM$, is trianalytic, i. e.
analytic with respect to all induced complex structures
(\ref{_tange_cone_underly_Proposition_}). 
It was proven in \cite{_Verbitsky:Deforma_}
(see also \ref{_triana_comple_geo_Proposition_}),
that trianalytic subvarieties are {\bf totally geodesic},
i. e. all geodesics in such a subvariety remain geodesics in
the ambient manifold. Since $T_xM$ is flat, its totally
geodesic subvariety $Z_xM$ is a union of planes. Now 
from the local homogeneity of singularities of $M$ it follows
that $M$ looks locally as its tangent cone, i. e. as
a union of non-singular hypercomplex varieties. 
This finishes the proof of Desingularization Theorem.

\subsection{Contents}

The paper is organized as follows.

\begin{itemize}

\item The Introduction is independent from the rest of this paper.

\item In Section \ref{_real_ana_Section_}, we recall some standard
definitions and results from the theory of real analytic
spaces. We define an almost complex real analytic space and show
that the complex structure on a complex variety can be recovered
from the corresponding almost complex structure on the underlying real
analytic variety.

\item Section \ref{_hyperka_Section_} contains some well-known results 
and definitions from the hyperk\"ahler geometry. We define trianalytic
subvarieties of hyperk\"ahler manifolds and show that trianalytic 
subvarieties naturally apppear in the complex geometry of hyperk\"ahler
manifolds. 

\item In Section \ref{_hypercomple_Section_}, 
we define a hypercomplex variety.
Examples of hypercomplex varieties include trianalytic
subvarieties of hyperk\"ahler manifolds
(\ref{_trianalytic_Definition_}) and moduli spaces of certain 
kinds of stable bundles over hyperk\"ahler manifolds
(Subsection \ref{_hyperholomo_Subsection_}). 
We study the tangent cone $Z_xM\subset T_xM$ 
of a hypercomplex variety, and show that it is a union 
of linear subspaces of the Zariski tanhgent space $T_xM$.

\item Section \ref{_LHS_Section_} deals with spaces 
having locally homogeneous singularities (SLHS). 
Roughly speaking, these are analytic spaces $M$ 
for which every point $x\in M$ has a 
system of coordinates $z_1, ..., z_n$ such that the
corresponding epimorphism of formal completions
\[ \C [[z_1, ...,  z_n]] \arrow \hat\calo_x M \]
has a homogeneous kernel (\ref{_SLHS_Definition_},
\ref{_locally_homo_coord_Claim_}).

We show that a space has LHS if it is endowed with a
system of infinitesimal automorphisms acting as dilatation
on the tangent spaces (\ref{_homogeni_LHS_Proposition_}).
We show that every hypercomplex variety
is naturally equipped with such a system of 
automorphisms, thus proving that it is SLHS
(\ref{_hyperco_SLHS_Theorem_}).

\item In Section \ref{_desingu_Section_},
we prove the desingularization theorem for hypercomplex varieties:
a normalization of a hypercomplex variety is smooth
and hypercomplex (\ref{_desingu_Theorem_}). This result
is used to show that all induced complex structures on a hypercomplex
variety are integrable (\ref{_all_indu_comple_integra_Theorem_}).

\item In Section \ref{_twistors_Section_}, we define a twistor
space of a hypercomplex variety. We show how to reconstruct a 
hypercomplex variety from its twistor space. We axiomatize this
situation, giving two sets of conditions which are satisfied
for twistor spaces of hypercomplex varieties: we define
{\bf twistor spaces of hypercomplex type}
(\ref{_twi_hyperco_type_Definition_}) and {\bf twistor
spaces of Deligne-Simpson type} (\ref{_Deligne_Simpson_Definition_}).

\item In Section \ref{_Deli_Si_equi_hyperco_Section_},
we prove that these sets of conditions are equivalent:
a variety is a twistor space of hypercomplex type if and
only if it is a twistor space of Deligne--Simpson type.

\item In Section \ref{_hype_type_equi_hype_Section_},
we show how to construct a hypercomplex variety starting
from an arbitrary twistor space of hypercomplex (or, equivalently,
Deligne-Simpson) type. This proves that a functor associating
to each hypercomplex variety a twistor space of hypercomplex type 
is an equivalence of categories.

\item In Section \ref{_twi_applications_Section_},
we give some applications of the equivalence of categories
constructed in Section \ref{_hype_type_equi_hype_Section_}.
We define hypercomplex spaces (\ref{_hyperco_spaces_Definition_}), 
thus generalizing the definition of hypercomplex varieties 
to spaces with nilpotents. We show that a quotient
of a hypercomplex variety by an action of a finite group
$G$ cannot be hypercomplex, unless $G$ acts freely.
Finally, we give another proof that the space of stable bundles
over a compact hyperk\"ahler manifold is hypercomplex,
assuming its Chern classes are ``suitable''
(for the precise definition of suitability and
the full statement, see Subsection 
\ref{_hyperholomo_Subsection_}; see also 
\cite{_Verbitsky:Hyperholo_bundles_}).

\end{itemize}


\section{Real analytic varieties and spaces}
\label{_real_ana_Section_}


In this section, we follow \cite{_GMT_}.

\subsection{Real analytic varieties and spaces: reduction, differentials}

Let $I$ be an ideal sheaf in the ring of real analytic functions
in an open ball $B$ in $\R^n$. The set of common zeroes 
of $I$ is equipped with a structure of ringed space,
with $\calo(B)/I$ as the structure sheaf. We denote
this ringed space by $Spec(\calo(B)/I)$.

\hfill

\definition
By a {\bf weak real analytic space} we understand a
ringed space which is locally isomorphic to
$Spec(\calo(B)/I)$, for some ideal $I \subset \calo(B)$.
A {\bf real analytic space} is a weak real analytic space
for which the structure sheaf is coherent
(i. e., locally finitely generated and presentable).

\hfill

For every real analytic variety $X$, there is
a natural sheaf morphism of evaluation,
$\calo(X) \stackrel{ev} \arrow C(X)$, 
where $C(X)$ is the sheaf of real analytic functions
on $X$. 

\hfill

\definition
A {\bf real analytic variety} is a weak real analytic
space for which the natural sheaf morphism $\calo(X) \arrow C(X)$
is injective. 

\hfill

Let $(X, \calo(X))$ be a real analytic space
and $N(X)\subset \calo(X)$ be the  kernel of 
the natural sheaf morphism $\calo(X) \arrow C(X)$.
Clearly, the ringed space $(X, \calo(X)/N(X))$
is a real analytic variety. This variety is
called {\bf a reduction of $X$}, denoted $X_r$.
The structure sheaf of $X_r$ is not necessarily coherent,
for examples see \cite{_GMT_}, III.2.15.

\hfill

For an ideal $I\subset \calo(B)$
we define the module of real analytic
differentials on $\calo(B)/I$ 
by 
\[ 
   \Omega^1(\calo(B)/I) = 
   \Omega^1(\calo(B))\bigg/\bigg ( I \cdot \Omega^1(\calo(B))+ d I \bigg),
\]
where $B$ is an open ball in $\R^n$, and
$\Omega^1(\calo(B))\cong \R^n \otimes \calo(B)$ is the module
of real analytic differentials on $B$. Patching this construction,
we define the sheaf of real analytic differentials on any real analytic
space. Likewise, one defines sheaves of
analytic differentials for complex varieties and in other
similar situations.

\subsection{Real analytic spaces underlying complex
analytic varieties}

Let $X$ be a complex analytic variety. The {\bf real analytic
space underlying $X$} (denoted by $X_\R$)
is the following object. By definition, $X_\R$ is a ringed space
with the same topology as $X$, but with a different structure
sheaf, denoted by $\calo_{X_\R}$. Let $i:\; U \hookrightarrow B^n$
be a closed complex analytic 
embedding of an open subset $U\subset X$ to 
an open ball $B^n\subset\C^n$,
and $I$ be an ideal defining $i(U)$. Then 
\[ \calo_{X_\R}\restrict U:=\; \calo_{B^n_\R}/Re(I)\]
is a quotient sheaf of the sheaf 
of real analytic functions on $B^n$ by the ideal $Re(I)$ generated
by the real parts of the functions $f\in I$.

Note that the real analytic
space underlying $X$ needs not be reduced, though it has no nilpotents
in the structure sheaf.
 
Consider the sheaf $\calo_X$ of holomorphic functions on $X$
as a subsheaf of the sheaf $C(X,\C)$ of continuous $\C$-valued
functions on $X$. The sheaf $C(X,\C)$ has a natural automorphism
$f\arrow \bar f$, where $\bar f$ is complex conjugation. 
By definition, the section $f$ of $C(X,\C)$  is called
{\bf antiholomorphic} if $\bar f$ is holomorphic.
Let $\calo_X$ be the sheaf of holomorphic functions,
and $\bar \calo_X$ be the sheaf of antiholomorphic
functions on $X$. Let $\calo_X \otimes_\C \bar\calo_X 
\stackrel i\arrow \calo_{X_\R}\otimes \C$ be the natural multiplication 
map. 

\hfill

\claim \label{_comple_real_ana_produ_Claim_}
Let $X$ be a complex variety, $X_\R$ the underlying
real analytic space.
Then the natural
 sheaf homomorphism $i:\; \calo_X \otimes_\C \bar\calo_X \arrow 
\calo_{X_\R}\otimes \C$
is injective. For each point $x\in X$, $i$ induces an isomorphism
on $x$-completions of $\calo_X \otimes_\C \bar\calo_X$
and $\calo_{X_\R}\otimes \C$.
 
{\bf Proof:} Clear from the definition. $\;\;\blacksquare$

\hfill

In the assumptions of \ref{_comple_real_ana_produ_Claim_},
let 
\[ \Omega^1(\calo_{X_\R}), \ \ 
   \Omega^1(\calo_X \otimes_\C \bar\calo_X),\ \  
   \Omega^1(\calo_{X_\R}\otimes \C)
\]
be the sheaves of real analytic
differentials associated with the corresponding
sheaves of rings. There is  a natural sheaf map

\begin{equation} \label{_Omega_X_R_and_Omega_X_Equation_}
  \Omega^1(\calo_{X_\R})\otimes \C = 
  \Omega^1(\calo_{X_\R}\otimes \C)\arrow 
   \Omega^1(\calo_X\otimes_\C \bar \calo_X),
\end{equation}
correspoding to the monomorphism
 
\[ 
   \calo_X\otimes_\C \bar \calo_X\hookrightarrow\calo_{X_\R}\otimes \C.
\]
 
\hfill

\claim \label{_differe_real_ana_and_co_ana_Claim_} 
Tensoring both sides of \eqref{_Omega_X_R_and_Omega_X_Equation_} 
by $\calo_{X_\R}\otimes \C$ produces an isomorphism
\[ \Omega^1(\calo_X\otimes_\C \bar \calo_X) 
   \bigotimes_{\calo_X\otimes_\C \bar \calo_X}\bigg(\calo_{X_\R}\otimes \C\bigg)   =\Omega^1(\calo_{X_\R}\otimes \C).
\]
 
{\bf Proof:} Clear. $\;\;\blacksquare$

\hfill

According to the general results about differentials
(see, for example, \cite{_Hartshorne:Alg_Geom_}, Chapter II,
Ex. 8.3), the sheaf
$\Omega^1(\calo_X\otimes_\C \bar \calo_X)$ admits a canonical
decomposition:
 
\[ \Omega^1(\calo_X\otimes_\C \bar \calo_X) =
   \Omega^1(\calo_X)\otimes_\C \bar \calo_X
   \oplus\calo_X\otimes_\C\Omega^1(\bar \calo_X).
\]
Let $\tilde I$ be an endomorphism of 
$\Omega^1(\calo_X\otimes_\C \bar \calo_X)$
which acts as a multiplication by $\1$ on
 
\[ \Omega^1(\calo_X)\otimes_\C \bar \calo_X
   \subset \Omega^1(\calo_X\otimes_\C \bar \calo_X)
\]
and as a multiplication by $-\1$ on
 
\[ \calo_X\otimes_\C\Omega^1(\bar \calo_X)
   \subset \Omega^1(\calo_X\otimes_\C \bar \calo_X).
\]
Let $\underline I$ be the corresponding 
$\calo_{X_\R}\otimes \C$-linear endomorphism of 
\[ \Omega^1(\calo_{X_\R})\otimes \C =
   \Omega^1(\calo_X\otimes_\C \bar \calo_X) 
   \otimes_{\calo_X\otimes_\C \bar \calo_X}
   \bigg(\calo_{X_\R}\otimes \C\bigg).
\]
A quick check shows that $\underline I$
is {\it real}, that is, comes from the 
$\calo_{X_\R}$-linear endomorphism of $\Omega^1(\calo_{X_\R})$.
Denote this $\calo_{X_\R}$-linear endomorphism by
\[ 
   I:\; \Omega^1(\calo_{X_\R})\arrow \Omega^1(\calo_{X_\R}), 
\]
$I^2=-1$. The endomorphism $I$ is called {\bf the complex structure
operator on the underlying real analytic space}. 
 In the case when $X$ is smooth, $I$ coinsides with
the usual complex structure operator on the cotangent bundle.

\hfill

\definition 
Let $M$ be a weak real analytic space, and
\[ I:\; \Omega^1(\calo_M)\arrow\Omega^1(\calo_M) \]
be an endomorphism satisfying $I^2=-1$. Then
$I$ is called {\bf an almost complex structure
on $M$}.

\subsection{Real analytic varieties and 
almost complex structures}

Let $B$ be an open ball in $\C^n$, and $X\subset B$ a closed complex
subvariety defined by an ideal $I \subset \calo_B$. Let 
$X_\R\subset B_\R$ be the underlying real analytic space,
and $X_\R^r\subset B_\R$ the underlying real analytic variety,
with the respective ideal sheaves denoted by
$I_\R$, $I^r_\R$. Consider the ideal
$I_\R^r \otimes \C \subset \calo_{B_\R}\otimes\C$.

\hfill

\lemma\label{_I_R^r_generators_Lemma_}
In the above assumptions, the ideal $I_\R^r\otimes \C$
is generated by elements
$f\cdot \bar g$, where $f, g\in \calo_B$ are 
holomorphic functions on $B$ satisfying $fg \in I$

{\bf Proof:} Clear. \endproof

\hfill

Let $F$ be a finite generated sheaf over a real analytic variety $Z$.
For all $z\in Z$, let ${\goth m}_z\subset \calo_Z$ be the ideal of
all functions vanishing in $z$. The vector space 
$ \cdot F / {\goth m}_z \cdot F$ is called
{\bf the fiber of $F$ in $z$}, denoted by $F\restrict z$.
Likewise we define the fiber $f\restrict z\in F\restrict z$
for a section $f$ of $F$.

\hfill

\lemma\label{_all_fibers_zero=>Fzero_Lemma_}
Let $F$ be a finite generated $\calo_Z$-sheaf 
over a real analytic variety $Z$, and $f$ its section.
Assume that $f\restrict z=0$ for all $z\in Z$. Then $f=0$.

\hfill

{\bf Proof:} Going to a closed 
subvariety if necessary, we may assume that
$Sup(f)=Z$. Let $Z_0\subset Z$ be an open subset such that
$F\restrict {Z_0}$ is free. Clearly, it suffices to show that
when $f\restrict{Z_0}$ is zero.
Since $F\restrict{Z_0}$ is free, $f\restrict{Z_0}$ is an $n$-tuple
of functions $(f_1, ... f_n)$, and $f\restrict z=0$
if and only if all $f_i$ take value $0$ at $z$. Applying
the definition of real analytic variety, we obtain $f=0$.
\endproof

\hfill

\proposition \label{_redu_differe_Lemma_}
Let $X_\R$ be a real analytic space underlying
a complex variety $X$, $X_\R^r$ be its reduction,
and $\Omega^1(X)$, $\Omega^1(X_\R^r)$ the corresponding sheaves
of real analytic differentials. Consider the natural map
\begin{equation} \label{_diffe_multipli_redu_arrow_Equation_}
    \Omega^1(X_\R)\otimes_{\calo_X}\calo_{X_\R^r}
    \stackrel{\phi_r} \arrow 
    \Omega^1(X_\R^r).
\end{equation}
Then $\phi_r$ is an isomorphism.

\hfill

{\bf Proof:} We work in notation introduced earlier in this section.
Consider the closed embedding $X_\R^r \hookrightarrow X_\R$. 
Let $N\subset \calo_{X_\R}$ be the ideal defining $X_\R^r$,
\[ N = \ker \bigg(\calo_{X_\R} \stackrel{ev} \arrow C(X) \bigg).\]
Clearly from the definitions,
\[ \Omega^1 (X_\R^r) = \Omega^1(X_\R) \otimes \calo_{X_\R^r}
   \bigg/ d N \otimes\calo_{X_\R^r}.
\]
To show that \eqref{_diffe_multipli_redu_arrow_Equation_}
is an isomorphism, it suffices to prove that the subsheaf
$dN \subset \Omega^1(X_\R)$, tensored by $\calo_{X_\R^r}$,
gives zero. By \ref{_all_fibers_zero=>Fzero_Lemma_},
it suffices to show that every section of 
\[ 
  d N \otimes\calo_{X_\R^r}\subset 
  \Omega^1(X_\R) \otimes \calo_{X_\R^r},
\] 
has zero fibers in $x$, for all points $x\in X$.

The fiber of $\Omega^1 X_\R$ in $x\in X$ is
${\goth m}_x/{\goth m}_x^2$, where 
${\goth m}_x\subset \calo_{X_\R}$ is the ideal generated
by all functions vanishing in $x$. For all $f\in \calo_{X_\R}$
such that $f(x)=0$, the fiber $df\restrict x$ of
$df\in\Omega^1 X_\R$ in $x$ coinsides with the class
of $f$ in ${\goth m}_x/{\goth m}_x^2$.

By \ref{_I_R^r_generators_Lemma_}, $N\otimes \C$ is 
generated by $f\bar g$,
$f, g \in \calo_B$, $fg\in I$. Therefore, 
$N\subset {\goth m}_x^2\otimes \C$ and by the above, the fiber 
$dN\restrict x$ iz zero.
This proves \ref{_redu_differe_Lemma_}. \endproof

\hfill

{}From \ref{_redu_differe_Lemma_}, it follows that
a a real analytic variety underlying a given complex variety
is equipped with a natural almost complex structure.
The corresponding operator is called {\bf the complex structure
operator in the underlying real analytic variety.}


\subsection{Integrability of almost complex structures}


\definition\label{_commu_w_comple_str_Definition_} 
Let $X$, $Y$ be complex analytic varieties, and 
\[ f:\; X_\R\arrow Y_\R\] be a morphism of underlying real
analytic varieties. Let 
$f^* \Omega^1_{Y_\R} \stackrel P\arrow \Omega^1_{X_\R}$ be the
natural map of sheaves of differentials associated with $f$.
Let 
 
\[ I_X:\; \Omega^1_{X_\R}\arrow \Omega^1_{X_\R}, \;\;\;
   I_Y:\; \Omega^1_{Y_\R}\arrow \Omega^1_{Y_\R} 
\]
be the complex structure operators, and
\[ f^* I_Y:\; f^*\Omega^1_{Y_\R}\arrow f^*\Omega^1_{Y_\R} \]
be $\calo_{X_\R}$-linear automorphism of 
$f^*\Omega^1_{Y_\R}$ defined as a pullback of $I_Y$.
We say that $f$ {\bf commutes with the complex structure}
if 
 
\begin{equation}\label{_commu_w_comle_Equation_}
   P\circ f^* I_Y = I_X \circ P.
\end{equation}
 \hfill

\theorem \label{_commu_w_comple_str_Theorem_} 
Let $X$, $Y$ be complex analytic varieties, and 
\[ f_\R:\; X_\R\arrow Y_\R\] be a morphism of underlying real
analytic varieties which commutes with the complex structure.
Then there exist a morphism $f:\; X\arrow Y$ of 
complex analytic varieties, such that $f_\R$ 
is its underlying morphism.

\hfill

{\bf Proof:} By Corollary 9.4, \cite{_Verbitsky:Deforma_}, the map
$f$, defined on the sets of points of $X$ and $Y$,
is meromorphic; to prove \ref{_commu_w_comple_str_Theorem_},
we need to show it is holomorphic. Let $\Gamma \subset X \times Y$
be the graph of $f$. Since $f$ is meromorphic, $\Gamma$ is 
a complex subvariety of $X\times Y$.
It will suffice to show that the natural projections
$\pi_1:\; \Gamma \arrow X$, $\pi_2:\; \Gamma \arrow Y$ are
isomorphisms. By \cite{_Verbitsky:Deforma_}, Lemma 9.12, 
the morphisms $\pi_i$ are flat. Since $f_\R$ induces isomorphism
of Zariski tangent spaces, same is true of $\pi_i$. Thus,
$\pi_i$ are unramified. Therefore, the maps $\pi_i$ are
etale. Since they are one-to-one on
points, $\pi_i$ etale implies $\pi_i$ is an isomorphism.
 $\;\;\blacksquare$
 
\hfill

\hfill

\definition \label{_integra_almost_comple_Definition_}
Let $M$ be a real analytic variety, and
\[ I:\; \Omega^1(\calo_M)\arrow\Omega^1(\calo_M) \]
be an endomorphism satisfying $I^2=-1$. Then
$I$ is called {\bf an almost complex structure
on $M$}. If there exist a structure $\goth C$ of complex variety
on $M$ such that $I$ appears as the complex structure operator
associated with $\goth C$, we say that $I$ is {\bf integrable}. 
\ref{_commu_w_comple_str_Theorem_} implies
that this complex structure is unique if it
exists.


\section{Hyperk\"ahler manifolds}
\label{_hyperka_Section_}


\subsection{Definitions}

This subsection contains a compression of 
the basic definitions from hyperk\"ahler geometry, found, for instance, in
\cite{_Besse:Einst_Manifo_} or in \cite{_Beauville_}.
 
\hfill
 
\definition \label{_hyperkahler_manifold_Definition_} 
(\cite{_Besse:Einst_Manifo_}) A {\bf hyperk\"ahler manifold} is a
Riemannian manifold $M$ endowed with three complex structures $I$, $J$
and $K$, such that the following holds.
 
\begin{description}
\item[(i)]  the metric on $M$ is K\"ahler with respect to these complex 
structures and
 
\item[(ii)] $I$, $J$ and $K$, considered as  endomorphisms
of a real tangent bundle, satisfy the relation 
$I\circ J=-J\circ I = K$.
\end{description}
 
\hfill 
 
The notion of a hyperk\"ahler manifold was 
introduced by E. Calabi (\cite{_Calabi_}).

\hfill
 
Clearly, a hyperk\"ahler manifold has the natural action of
the quaternion algebra ${\Bbb H}$ on its real tangent bundle $TM$. 
Therefore its complex dimension is even.
For each quaternion $L\in \Bbb H$, $L^2=-1$,
the corresponding automorphism of $TM$ is an almost complex
structure. It is easy to check that this almost 
complex structure is integrable (\cite{_Besse:Einst_Manifo_}).
 
\hfill
 
\definition \label{_indu_comple_str_Definition_} 
Let $M$ be a hyperk\"ahler manifold, $L$ a quaternion satisfying
$L^2=-1$. The corresponding complex structure on $M$ is called
{\bf an induced complex structure}. The $M$ considered as a complex
manifold is denoted by $(M, L)$.
 
\hfill
 
Let $M$ be a hyperk\"ahler manifold. We identify the group $SU(2)$
with the group of unitary quaternions. This gives a canonical 
action of $SU(2)$ on the tangent bundle, and all its tensor
powers. In particular, we obtain a natural action of $SU(2)$
on the bundle of differential forms. 

\hfill

\lemma \label{_SU(2)_commu_Laplace_Lemma_}
The action of $SU(2)$ on differential forms commutes
with the Laplacian.
 
{\bf Proof:} This is Proposition 1.1
of \cite{_Verbitsky:Hyperholo_bundles_}. \endproof
 
Thus, for compact $M$, we may speak of the natural action of
$SU(2)$ in cohomology.

 
\subsection{Trianalytic subvarieties in compact hyperk\"ahler
manifolds.}
 
 
In this subsection, we give a definition and a few basic properties
of trianalytic subvarieties of hyperk\"ahler manifolds. 
We follow \cite{_Verbitsky:Symplectic_II_}.
 
\hfill
 
Let $M$ be a compact hyperk\"ahler manifold, $\dim_\R M =2m$.
 
\hfill
 
\definition\label{_trianalytic_Definition_} 
Let $N\subset M$ be a closed subset of $M$. Then $N$ is
called {\bf trianalytic} if $N$ is a complex analytic subset 
of $(M,L)$ for any induced complex structure $L$.
 
\hfill
 
Let $I$ be an induced complex structure on $M$,
and $N\subset(M,I)$ be
a closed analytic subvariety of $(M,I)$, $dim_\C N= n$.
Denote by $[N]\in H_{2n}(M)$ the homology class 
represented by $N$. Let $\inangles N\in H^{2m-2n}(M)$ denote 
the Poincare dual cohomology class. Recall that
the hyperk\"ahler structure induces the action of 
the group $SU(2)$ on the space $H^{2m-2n}(M)$.
 
\hfill
 
\theorem\label{_G_M_invariant_implies_trianalytic_Theorem_} 
Assume that $\inangles N\in  H^{2m-2n}(M)$ is invariant with respect
to the action of $SU(2)$ on $H^{2m-2n}(M)$. Then $N$ is trianalytic.
 
{\bf Proof:} This is Theorem 4.1 of 
\cite{_Verbitsky:Symplectic_II_}.
\endproof
 
\hfill

\remark \label{_triana_dim_div_4_Remark_}
Trianalytic subvarieties have an action of quaternion algebra in
the tangent bundle. In particular,
the real dimension of such subvarieties is divisible by 4.

\hfill

Let $M$ be a hyperk\"ahler manifold, $\c R$ the set 
of induced complex structures. The following
result implies that for generic
$I\in \c R$, all complex subvarieties of $(M, I)$ 
are trianalytic in $M$.

\hfill

\definition
Let $M$ be a compact hyperk\"ahler manifold, $\c R$ the set 
of induced complex structures. The complex structure
$I\subset \c R$ is called {\bf of general type with respect
to the hyperk\"ahler structure} if
all rational $(p,p)$-classes 
\[ \omega \in 
  \bigoplus\limits_p H^{p,p}(M)\cap H^{2p}(M,\Z)\subset H^*(M)
\] 
are $SU(2)$-invariant.

\hfill

\proposition
Let $M$ be a compact hyperk\"ahler manifold, $\c R$ the set 
of induced complex structures. Then for all $I\in \c R$
except a countable subset, $I$ is of general type.

{\bf Proof:} This is Proposition 2.2 from
\cite{_Verbitsky:Symplectic_II_}. \endproof


\subsection{Totally geodesic submanifolds.}


\nopagebreak
\hspace{5mm}
\proposition \label{_comple_geodesi_basi_Proposition_}
Let $X \stackrel \phi\hookrightarrow M$ be an embedding of Riemannian 
manifolds (not necessarily compact) compatible with the Riemannian
structure.
 Then the following conditions are equivalent.
 
\begin{description}
\item[(i)] Every geodesic line in $X$ is geodesic in $M$.
 
\item[(ii)] Consider the Levi-Civita connection $\nabla$ on $TM$,
and restriction of $\nabla$ to $TM \restrict{X}$. Consider the
orthogonal decomposition 
\begin{equation} \label{TM_decompo_Equation_} 
   TM\restrict{X} = TX \oplus TX^\bot. 
\end{equation}
Then, this decomposition is preserved by the connection $\nabla$.
\end{description}
 
{\bf Proof:} Well known; see, for instance, 
\cite{_Besse:Einst_Manifo_}.
\blacksquare
 
\hfill

\proposition \label{_triana_comple_geo_Proposition_} 
Let $X\subset M$ be a trianalytic submanifold of a hyperk\"ahler
manifold $M$, where $M$ is not necessarily compact. Then
$X$ is totally geodesic.

{\bf Proof:} This is \cite{_Verbitsky:Deforma_}, Corollary 5.4.
\endproof


\section{Hypercomplex varieties}
\label{_hypercomple_Section_}


\subsection{Definition and examples}

\definition\label{_almost_hyperco_Definition_}
Let $M$ be a real analytic variety equipped with almost
complex structures $I$, $J$ and $K$, such that
$I\circ J = -J \circ I = K$. Then $M$ is called 
{\bf an almost hypercomplex variety.} 

\hfill

An almost hypercomplex variety is equipped with an action of
quaternion algebra in its differential sheaf. Each quaternion
$L\in \Bbb H$, $L^2=-1$ defines an almost complex structure
on $M$. Such an almost complex structure is called
{\bf induced by the hypercomplex structure}.

\hfill

\definition\label{_hyperco_Definition_}
Let $M$ be an almost hypercomplex variety. We say that
$M$ is {\bf hypercomplex} if there exist a pair of induced complex
structures $I_1, I_2\in \Bbb H$, $I_1\neq \pm I_2$, such that
$I_1$ and $I_2$ are integrable.

\hfill

{\bf Caution:} Not everything which looks hypercomplex
satisfies the conditions of \ref{_hyperco_Definition_}.
Take a quotient $M/G$ of a hypercomplex manifold by an action 
of a finite group $G$, acting compatible with the hyperk\"ahler
structure. Then $M/G$ is {\it not} hypercomplex, unless
$G$ acts freely (\ref{_quotie_not_hyperco_Proposition_}).

\hfill

\noindent
\claim \label{_hyperka_hyperco_Claim_} 
Let $M$ be a hyperk\"ahler manifold. Then $M$ is hypercomplex.
{\bf Proof:} Let $I$, $J$ be induced complex structures.
We need to identify $(M, I)_\R$ and $(M,J)_\R$ in a natural way.
These varieties are canonically identified as $C^\infty$-manifolds;
we need only to show that this identification is real analytic.
This is \cite{_Verbitsky:Deforma_}, Proposition 6.5. \endproof

\hfill

\remark\label{_triana_hyperco_Remark_}
Trianalytic subvarieties of hyperk\"ahler manifolds
are obviously hypercomplex. 
Define trianalytic subvarieties of hypercomplex varieties
as subvarieties which are complex analytic with respect to 
all integrable induced
complex structures. Clearly, trianalytic subvarieties of
 hypercomplex varieties are equipped with a natural hypercomplex
structure. Another example of a hypercomplex
variety is given in Subsection \ref{_hyperholomo_Subsection_}.
For additional examples, see \cite{_Verbitsky:Deforma_}.

\subsection{Tangent cone of a hypercomplex variety.}

Let $M$ be a hypercomplex variety, $I$ an integrable 
induced complex structure, $\tilde Z_x(M,I)$ be 
the Zariski tangent cone to $(M,I)$ in $x\in M$
and $Z_x(M,I)$ its reduction. 
Consider $Z_x(M,I)$ as a closed subvariety in 
the Zariski tangent space $T_xM$. The space $T_xM$  has a natural
quaternionic structure and admits a
compatible metric. This makes $T_xM$ into a 
hyperk\"ahler manifold, isomorphic to ${\Bbb H}^n$. 

\hfill

\noindent\theorem \label{_cone_hype_Theorem_}
Under these assumptions, the following assertions hold:
\begin{description}
\item[(i)] The subvariety $Z_x(M,I)\subset T_x M$ is independent
from the choice of integrable induced complex structure $I$.
\item [(ii)] Moreover, $Z_x(M,I)$ is a trianalytic subvariety
of $T_x M$.
\end{description}

{\bf Proof:} Clearly,
\ref{_cone_hype_Theorem_} (ii) follows from
\ref{_cone_hype_Theorem_} (i).
\ref{_cone_hype_Theorem_} (i) is 
directly implied by the following general result.

\hfill

\proposition\label{_tange_cone_underly_Proposition_} 
Let $M$ be a complex variety, $x\in X$ a point, and $Z_xM\subset T_xM$
be the reduction of the Zariski tangent cone to $M$ in $x$, considered
as a closed subvariety of the Zariski tangent space $T_xM$.
Let $Z_x M_\R \subset T_x M_\R$ be the Zariski tangent cone for the
underlying real analytic space $M_\R$. Then 
$(Z_x M)_\R \subset (T_x M)_\R = T_x M_\R$ coinsides with $Z_x M_\R$.

{\bf Proof:} For each $v\in T_x M$, the point $v$ belongs to
$Z_x M$ if and only if there exist a real analytic path
$\gamma:\; [0, 1] \arrow M$, $\gamma(0)=x$ satisfying 
$\frac{d\gamma}{dt}=v$. The same holds true for $Z_x M_\R$.
Thus, $v\in Z_x M$ if and only if $v\in Z_x M_\R$. \endproof

\hfill

The following theorem shows that the Zariski tangent cone
$Z_xM\subset T_x M$ is a union of planes $L_i\subset T_x M$.

\hfill

\theorem \label{_cone_flat_Theorem_}
Let $M$ be a hypercomplex variety, $I$ an induced complex
structure and $x\in M$ a point. 
Consider the reduction of the
Zariski tangent cone (denoted by $Z_x M$) as a subvariety of the
quaternionic space $T_x M$. Let $Z_x(M, I)= \cup L_i$ 
be the irreducible decomposition of the complex variety $Z_x(M,I)$.
Then

\begin{description}
\item[(i)] The decomposition $Z_x(M, I)= \cup L_i$ 
is independent from the choice of induced complex structure $I$.
\item[(ii)] For every $i$, the variety 
$L_i$ is a linear subspace of $T_x M$,
invariant under quaternion action.
\end{description}

{\bf Proof:}  Let $L_i$ be an irreducible component of
$Z_x(M, I)$, $Z_x^{ns}(M,I)$ be the non-singular part
of $Z_x(M,I)$, and $L_i^{ns}:=Z_x^{ns}(M,I) \cap L_i$.
Then $L_i$ is a closure of $L_i^{ns}$ in $T_xM$.
Clearly from \ref{_cone_hype_Theorem_}, $L_i^{ns}(M)$ is a hyperk\"ahler
submanifold in $T_xM$. By 
\ref{_triana_comple_geo_Proposition_}, $L_i^{ns}$ is totally
geodesic. A totally geodesic submanifold of a flat manifold is
again flat. Therefore, $L_i^{ns}$ is an open subset of a linear
subspace $\tilde L_i\subset T_xM$. Since $L_i^{ns}$ is a hyperk\"ahler
submanifold, $\tilde L_i$ is invariant with respect to quaternions.
The closure $L_i$ of $L_i^{ns}$ is a complex analytic subvariety
of $T_x(M,I)$. Therefore, $\tilde L_i = L_i$. This proves 
\ref{_cone_flat_Theorem_} (ii). From the above argument, it is
clear that $Z_x^{ns}(M,I)= \coprod L_i^{ns}$ (disconnected sum).
Taking connected components of $Z_x^{ns}M$ for each induced 
complex structure, we obtain the same decomposition
$Z_x(M, I)= \cup L_i$, with $L_i$ being closured of connected components.
This proves \ref{_cone_flat_Theorem_} (ii). \endproof


\section[Hypercomplex varieties have locally homogeneous singularities]{Hypercomplex varieties \\have locally homogeneous singularities}
\label{_LHS_Section_}


This section follows \cite{_Verbitsky:DesinguII_}.

\subsection{Spaces with locally homogeneous singularities.}

\noindent
\definition
(local rings with LHS)
Let $A$ be a local ring. Denote by $\goth m$ its maximal ideal.
Let $A_{gr}$ be the corresponding associated graded ring
for the $\goth m$-adic filtration.
Let $\hat A$, $\widehat{A_{gr}}$ be the $\goth m$-adic completion
of $A$, $A_{gr}$. Let $(\hat A)_{gr}$, $(\widehat{A_{gr}})_{gr}$ 
be the associated graded rings, which are naturally isomorphic to
$A_{gr}$. We say that $A$ {\bf has locally homogeneous singularities}
(LHS)
if there exists an isomorphism $\rho:\; \hat A \arrow \widehat{A_{gr}}$
which induces the standard isomorphism 
$i:\; (\hat A)_{gr}\arrow (\widehat{A_{gr}})_{gr}$ on associated
graded rings.

\hfill

\definition\label{_SLHS_Definition_}
(SLHS)
Let $X$ be a complex or real analytic space. Then 
$X$ is called {\bf a space with locally homogeneous singularities}
(SLHS) if for each $x\in X$, the local ring $\calo_x X$ 
has locally homogeneous singularities.

\hfill
 
The following claim might shed a light on the origin of the term
``locally homogeneous singularities''.

\hfill

\claim \label{_locally_homo_coord_Claim_}
Let $A$ be a complete local Noetherian ring over $\C$,
with a residual field $\C$. 
Then the following statements are equivalent
\begin{description}
\item[(i)] $A$ has locally homogeneous singularities
\item[(ii)] There exist
a surjective
ring homomorphism $\rho:\; \C[[x_1, ... , x_n]] \arrow A$, 
where $\C[[x_1, ... , x_n]]$ is the ring of power series,   
and the ideal $\ker \rho$ is homogeneous in $\C[[x_1, ... , x_n]]$.
\end{description}

{\bf Proof:} Clear. \endproof

\hfill

\definition
Let $M$ be a hypercomplex variety. Then $M$ is called 
a space with locally homogeneous singularities (SLHS) if 
 for all integrable induced complex structures
$I$ the $(M, I)$ is SLHS.

\subsection{Complete rings with automorphisms}

\definition \label{_homogeni_automo_Definition_}
Let $A$ be a local Noetherian ring over $\C$,
with a residual field $\C$, equipped with an
automorphism
$e:\; A \arrow A$. Let $\goth m$ be a maximal ideal of $A$.
Assume that $e$ acts on $\goth m /\goth m^2$ as a multiplication 
by $\lambda\in \C$, $|\lambda|< 1$. Then $e$ is called {\bf a
homogenizing automorphism of $A$}.

\hfill

The aim of the present subsection is to prove the following statement.

\hfill

\proposition \label{_homogeni_LHS_Proposition_}
Let $A$ be a complete Noetherian ring over $\C$,
with a residual field $\C$, equipped with a
homogenizing authomorphism $e:\; A \arrow A$. Then there exist
a surjective 
ring homomorphism $\rho:\; \C[[x_1, ... , x_n]] \arrow A$, 
such that the ideal
$\ker \rho$ is homogeneous in $\C[[x_1, ... , x_n]]$.
In particular, $A$ has locally homogeneous 
singularities.\footnote{See \ref{_locally_homo_coord_Claim_}
for  LHS property in terms of coordinate systems.}

\hfill

This statement is well known. A reader who knows its proof
should skip the rest of this section. 

\hfill

\proposition \label{_homogeni_auto_then_basis_Proposition_}
Let $A$ be a complete Noetherian ring over $\C$,
with a residual field $\C$, equipped with a
homogenizing authomorphism $e:\; A \arrow A$. Then there exist
a system of ring elements
\[ 
    f_1 , ..., f_m \in \goth m, \ \ m = \dim_\C\goth m /\goth m^2,
\]
which generate $\goth m /\goth m^2$, and such that $e(f_i) = \lambda f_i$.

\hfill

{\bf Proof:}
Let $\underline a\in\goth m /\goth m^2$.
Let $a\in \goth m$ be a representative of $\underline a$ in $\goth m$.
To prove \ref{_homogeni_auto_then_basis_Proposition_}
it suffices to find $c \in \goth m^2$, such that 
$e(a-c) = \lambda a -\lambda c$. Thus, we need to solve an equation
\begin{equation}\label{_a_through_a_Equation_}
 \lambda c - e(c) = e(a) - \lambda(a). 
\end{equation}
Let $r:= e(a)-\lambda a$. Clearly, $r\in \goth m ^2$.
A solution of \eqref{_a_through_a_Equation_}
is provided by the following lemma.

\hfill

\lemma \label{_e-lambda_invertible_Lemma_}
In assumptions of \ref{_homogeni_auto_then_basis_Proposition_},
let $r\in \goth m^2$. Then, the equation
\begin{equation}\label{_finding_eigen_Equation_}
e(c) - \lambda c = r
\end{equation}
has a unique solution $c \in \goth m^2$.

\hfill

{\bf Proof:} We need to show that the operator
$P:= (e-\lambda)\restrict{\goth m^2}$
is invertible. Consider the $\goth m$-adic filtration 
$\goth m^2 \supset \goth m^3 \supset ...$ on $\goth m^2$.
Clearly, $P$ preserves this filtration. Since $\goth m^2$
is complete with respect to the adic filtration,
it suffices to show that $P$ is invertible on the
successive quotients. The quotient $\goth m^2/\goth m^i$ is 
finite-dimensional, so to show that $P$ is invertible it suffices
to calculate the eigenvalues. Since $e$ is an automorphism,
restriction of $e$ to $\goth m^i/\goth m^{i-1}$ is a multiplication
by $\lambda^i$. Thus, the eigenvalues of $e$ on $\goth m^2/\goth m^i$ 
range from $\lambda^2$ to $\lambda^{i-1}$. Since $|\lambda|>|\lambda|^2$,
all eigenvalues of $P\restrict{\goth m^2/\goth m^i}$ are 
non-zero and the restriction of $P$ to $\goth m^2/\goth m^i$ is invertible.
This proves \ref{_e-lambda_invertible_Lemma_}.
$\blacksquare$

\hfill

{\bf The proof of \ref{_homogeni_LHS_Proposition_}.}
Consider the map \[ \rho:\; \C[[x_1, ... x_m]] \arrow A,\ \ 
\rho(x_i) = f_i,\] where $f_1, ... , f_m$ is the system of functions
constructed in \ref{_homogeni_auto_then_basis_Proposition_}.
By Nakayama's lemma, $\rho$
is surjective. 

Let $e_\lambda:\; \C[[x_1, ... x_m]] \arrow\C[[x_1, ... x_m]] $ be the 
automorphism
mapping $x_i $ to $\lambda x_i$. Then, the diagram
\[\begin{CD} \C[[x_1, ... x_m]] @>{\rho}>> A \\
                @V{e_\lambda}VV @VV{e}V\\
                \C[[x_1, ... x_m]] @>{\rho}>> A 
\end{CD}
\]
is by construction commutative.
Therefore, the ideal $I= \ker \rho$ is preserved by $e_\lambda$.
Clearly, every $e_\lambda$-preserved 
ideal $I\subset \C[[x_1, ... x_m]]$ is homogeneous.
\ref{_homogeni_LHS_Proposition_} is proven. \endproof

\subsection{Automorphisms of local rings of hypercomplex varieties}

Let $M$ be a hypercomplex variety, $x\in M$ a point,
$I$ an integrable
induced complex structure. Let $A_I:= \hat \calo_x(M,I)$
be the adic completion of the local ring $\calo_x(M,I)$ of $x$-germs
of holomorphic functions on the complex variety $(M,I)$.
Clearly, the sheaf ring of the antiholomorphic functions on $(M,I)$
coinsides with $\calo_x(M,-I)$. Thus, the corresponding completion
ring is $A_{-I}$. 
The isomorphism of 
\ref{_comple_real_ana_produ_Claim_} produces a natural epimorphism
\begin{equation}\label{_co_ana_and_rea_isom_Equation_}
\widehat{A_I \otimes_\C A_{-I}} \arrow A_\R,
\end{equation}
where \[ A_\R := \widehat{\calo_x(M_\R)\otimes_\R \C}\] is the $x$-completion
of the ring of germs of real analytic complex-valued functions on $M$.
Consider the natural quotient map \[ p:\;A_{-I}\arrow \C.\]

Consider the natural epimorphism of complete rings
\begin{equation}\label{_epi_from_pro_to_A_I_Equation_}
   \widehat{A_I \otimes_\C A_{-I}} \arrow A_I,\ \  
   a\otimes b \mapsto a\otimes p(b),
\end{equation}
where $a\in A_I$, $b\in A_{-I}$, and
\[ a\otimes b\in{A_I \otimes_\C A_{-I}}.\]

\hfill

\lemma\label{_p_zero_on_kernel_of_multi_Lemma_}
The kernel of
\eqref{_epi_from_pro_to_A_I_Equation_} contains the kernel of
\eqref{_co_ana_and_rea_isom_Equation_}.

\hfill

{\bf Proof:} Consider an epimorphism $\phi:\; B_x \arrow A_I$
where \[ B_x = \C [[z_1, ... , z_n]].\] Let ${\goth I}\subset B_x$ be the 
kernel of $\phi$. 
By \ref{_I_R^r_generators_Lemma_}, the ring
$A_\R$ is naturally isomorphic to
\[ 
   (B_x)_\R = \C [[z_1, ... , z_n, \bar z_1, ... \bar z_n]]/ {\goth I}_\R, 
\]
where ${\goth I}_\R$ is an ideal generated by all the products
\[ f(z_1, ... z_n)\cdot\bar g(\bar z_z , ... \bar z_n),\]
such that $fg \in \goth I$. Likewise,
$\widehat{A_I \otimes_\C A_{-I}}$ is a quotient of
\[ (B_x)_\R = \C [[z_1, ... , z_n, \bar z_1, ... \bar z_n]] \]
by the ideal 
\[ {\goth I} \cdot \C [[\bar z_1, ... \bar z_n]] +
   \C [[z_1, ... , z_n]\cdot  \bar {\goth I}.
\]
Slightly abusing the notation, we denote 
the corresponding quotient map by
\[ \phi:\; \C [[z_1, ... , z_n, \bar z_1, ... \bar z_n]]
   \arrow \widehat{A_I \otimes_\C A_{-I}}.
\]

Let $a\in \widehat{A_I \otimes_\C A_{-I}}$
be an element which is mapped to zero by the map 
\eqref{_co_ana_and_rea_isom_Equation_}. Then 
$a$ is a linear combination of 
$\phi(f_i\bar g_i)$, for $f_i, g_i\in \C [[z_1, ... , z_n]]$,
$fg\in \goth I$. Therefore, it suffices to show that $a$
lies in the kernel of \eqref{_epi_from_pro_to_A_I_Equation_}
for $a = \phi(f \bar g)$.
Either $g$ is invertible and $f\in \goth I$, or 
$g(0, ... 0)=0$. In the first case,
$f\bar g\in {\goth I}\otimes \C [[\bar z_1, ... \bar z_n]]$, so 
$\phi(f\bar g)=0$. In the second case, $p(\bar g) =0$,
so $1\otimes p(a)=0$ and
$a$ lies in the kernel of the map
\eqref{_epi_from_pro_to_A_I_Equation_}.
This proves \ref{_p_zero_on_kernel_of_multi_Lemma_}.
\endproof

\hfill

Consider the diagram 
\[\begin{CD}
\widehat{A_I \otimes_\C A_{-I}} @>>> A_I\\
\searrow \\
\ \ \ \  \ \ \ \  \  A_\R
\end{CD}
\]
formed from the arrows of
\eqref{_epi_from_pro_to_A_I_Equation_} and 
\eqref{_co_ana_and_rea_isom_Equation_}.
By \ref{_p_zero_on_kernel_of_multi_Lemma_}, there exists
an epimorphism $e_I:\; A_\R \arrow A_I$, making this diagram
commutative. Let
$i_I:\; A_I \hookrightarrow A_\R$ be the natural embedding 
\[ a \mapsto a\otimes 1\in\widehat{A_I \otimes_\C A_{-I}}.\]
For an integrable induced complex structure $J$, we define
$A_J$, $A_{-J}$, $i_J$, $e_J$ likewise.

Let $\Psi_{I,J}:\; A_I \arrow A_I$ be the composition

\[ A_I \stackrel {i_I}\arrow A_\R\stackrel {e_J}\arrow A_J
   \stackrel {i_J}\arrow A_\R\stackrel {e_I}\arrow A_I.
\]
Clearly, for $I=J$, the ring morphism $\Psi_{I,J}$ is identity,
and for $I=-J$, $\Psi_{I,J}$ is an augmentation map.

\hfill

\proposition \label{_homogenizing_Proposition_}
Let $M$ be a hypercomplex variety, $x\in M$ a point, and
$I$, $J$ induced complex structures, such that $I\neq J$ and
$I\neq -J$. Consider the map  $\Psi_{I,J}:\; A_I \arrow A_I$
defined as above. Then $\Psi_{I,J}$ is a homogenizing
automorphism of $A_I$.\footnote{For the definition
of a homogenizing automorphism, see \ref{_homogeni_automo_Definition_}.}

\hfill

{\bf Proof:}
Let $d\Psi$ be the differential of $\Psi_{I,J}$, that is, 
the restriction of $\Psi_{I,J}$ to $*\goth m/\goth m^2)^*$,
where $\goth m$ is the maximal ideal of $A_I$. 
By Nakayama's lemma, to prove that
$\Psi_{I,J}$ is an automorphism it suffices to show that $d\Psi$
is invertible. To prove  that $\Psi_{I,J}$ is homogenizing, 
we have to show that $d\Psi$ is a multiplication by a complex
number $\lambda$, $|\lambda|<1$. As usually, we denote the real analytic
variety underlying $M$ by $M_\R$. 
Let $T_I$, $T_J$, $\underline {T}_\R$ be the Zariski 
tangent spaces to $(M,I)$, $(M,J)$ and $M_\R$, respectively, 
in $x\in M$. Consider the complexification
  $T_\R:= \underline {T}_\R\otimes \C$, which is a 
Zariski tangent space to the local ring $A_\R$.
To compute $d\Psi:\; T_I \arrow T_I$, we need
to compute the differentials of $e_I$, $e_J$, $i_I$, $i_J$,
i. e., the restrictions of the homomorphisms
$e_I$, $e_J$, $i_I$, $i_J$ to the Zariski tangent spaces
$T_I$, $T_J$, $T_\R$.
Denote these differentials by $de_I$, $de_J$, $di_I$, $di_J$.

\hfill

\lemma \label{_i_e_through_Hodge_Lemma_}
Let $M$ be a hypercomplex variety, $M_\R$ the associated real analytic
variety, $x\in M$ a point. Consider the space $T_\R := T_x (M_\R)\otimes
\C$. For an induced complex structure $I$, consider the Hodge 
decomposition
$T_\R= T^{1,0}_I \oplus  T^{0,1}_I$. In our previous notation,
$T_I^{1,0}$ is $T_I$. Then, $di_I$ is the natural embedding
of $T_I = T_I^{1,0}$ to $T_\R$, and $de_I$ is the natural
projection of $T_\R= T^{1,0}_I \oplus  T^{0,1}_I$ to
$T_I^{1,0}=T_I$.

{\bf Proof:} Clear. \endproof

\hfill

We are able now to describe the map $d\Psi:\; T_I \arrow T_I$
in terms of the quaternion action. Recall that the space $T_I$
is equipped with a natural $\R$-linear
quaternionic action. For each quaternionic
linear space $\underline V$ and each quaternion $I$, $I^2=-1$, $I$ defines a
complex structure in $\underline V$. Such a complex structure is called
{\bf induced by the quaternionic structure}. 

\hfill

\lemma \label{_Psi_through_quate_Lemma_}
Let $\underline V$ be a space with quaternion action, and 
$V:= \underline V \otimes \C$ its complexification.
For each induced complex structure $I\in {\Bbb H}$,
consider the Hodge decomposition $V:= V_I^{1,0} \oplus V_I^{0,1}$.
For induced complex structures $I, J\in \Bbb H$,
let $\Phi_{I,J}(V)$ be the composition of the natural embeddings
and projections
\[ 
   V_I^{1,0} \arrow V \arrow  V_J^{1,0} \arrow V \arrow  V_I^{1,0}.
\]
Using the natural identification $\underline V \cong  V_I^{1,0}$,
we consider $\Phi_{I,J}(V)$ as an $\R$-linear
automorphism of the space $\underline V$.
Then, applying the operator $\Phi_{I,J}(V)$  to the 
quaternionic space $T_I$, we obtain the operator $d\Psi$
defined above.

{\bf Proof:} Follows from \ref{_i_e_through_Hodge_Lemma_}
\endproof

\hfill

As we have seen, to prove \ref{_homogenizing_Proposition_}
it suffices to show that $d\Psi$ is a multiplication by a non-zero
complex number $\lambda$, $|\lambda| < 1$.
Thus, the proof of \ref{_homogenizing_Proposition_} is finished with the
following lemma.

\hfill

\lemma\label{_compu_of_Psi_for_qua_Lemma_}
In assumptions of \ref{_Psi_through_quate_Lemma_}, consider the 
map \[ \Phi_{I,J}(V):\; V_I^{1,0} \arrow V_I^{1,0}.\] Then $\Phi_{I,J}(V)$
is a multiplication by a complex number $\lambda$. 
Moreover, $\lambda$ is a non-zero number unless $I=-J$,
and $|\lambda|< 1$ unless $I=J$.

\hfill

{\bf Proof:}
Let $\underline V= \oplus \underline V_i$ be a decomposition of $V$ into
a direct sum of $\Bbb H$-linear spaces. Then, the operator $\Phi_{I,J}(V)$
can also be decomposed: $\Phi_{I,J}(V) = \oplus \Phi_{I,J}(V_i)$.
Thus, to prove \ref{_compu_of_Psi_for_qua_Lemma_} it suffices
to assume that $\dim_{\Bbb H} \underline V=1$.
Therefore, we may identify $\underline V$ with the space
$\Bbb H$, equipped with the right action of quaternion
algebra on itself.

Consider the left action of $\Bbb H$ on $\underline V = \Bbb H$.
This action commutes with the right action of $\Bbb H$ on $\underline V$.
Consider the corresponding action 
\[ 
   \rho:\; SU(2) \arrow \End(\underline V)
\] of the group of unitary
quaternions ${\Bbb H}^{un}=SU(2)\subset \Bbb H$ on $\underline V$. 
Since $\rho$ commutes with the quaternion action, 
$\rho$ preserves $V^{1,0}_I \subset V$, for every
induced complex structure $I$. In the same way, for each $g\in SU(2)$,
the endomorphism $\rho(g)\in \End(V^{1,0}_I)$
commutes with $\Phi_{I,J}(V)$.

Consider the 2-dimensional $\C$-vector space $V^{1,0}_I$ 
as a representation of $SU(2)$. Clearly, $V^{1,0}_I$ 
is an irreducible representation. Thus, by Schur's lemma,
the automorphism $\Phi_{I,J}(V)\in \End(V^{1,0}_I))$ 
is a multiplication by a complex constant $\lambda$.
The bounds $0< |\lambda| < 1$ are implied by the 
following elementary argument. The composition
$i_I \circ e_J$ applied to a vector 
$v\in V_I^{1,0}$ is a projection of $v$ to $V_J^{1,0}$
along $V_J^{0,1}$. Consider the natural Euclidean metric
on $V = \Bbb H$. Clearly, the decomposition 
$V = V_J^{1,0}\oplus V_J^{0,1}$ is orthogonal.
Thus, the composition $i_I \circ e_J$ is an orthogonal
projection of $v\in V_I^{1,0}$ to $V_J^{1,0}$.
Similarly, the composition $i_J \circ e_I$ is an orthogonal
projection of $v\in V_J^{1,0}$ to $V_I^{1,0}$.
Thus, the map $\Phi_{I,J}(V)$ is an orthogonal projection 
from $V_I^{1,0}$ to $V_J^{1,0}$ and back to $V_I^{1,0}$.
Such a composition always decreases a length of vectors,
unless $V_I^{1,0}$ coincides with $V_J^{1,0}$, in which
case $I=J$. Also, unless $V_I^{1,0} =  V_J^{0,1}$,
$\Phi_{I,J}(V)$ is non-zero; in the later case, $I = -J$.
\ref{_homogenizing_Proposition_}
is proven. \endproof

\hfill

{}From \ref{_homogenizing_Proposition_} and 
\ref{_homogeni_LHS_Proposition_}, we obtain the 
following theorem.

\hfill

\theorem\label{_hyperco_SLHS_Theorem_}
(the main result of this section)
Let $M$ be a hypercomplex variety. Then
$M$ is a space with locally homogeneous singularities
(SLHS).

\endproof


\section{Desingularization of hypercomplex varieties}
\label{_desingu_Section_}


\subsection{The proof of desingularization theorem}

\proposition\label{_normali_smooth_Corollary_}
Let $M$ be a hypercomplex 
variety, and $I$ an integrable induced 
complex structure.
Then the normalization of $(M,I)$ is smooth.

\hfill

{\bf Proof:} The normalization of $Z_xM$ is smooth by 
\ref{_cone_flat_Theorem_}. The 
normalization is compatible with the adic
completions (\cite{_Matsumura:Commu_Alge_},
Chapter 9, Proposition 24.E). Therefore, the integral closure
of the completion of $\calo_{Z_xM}$ is a regular ring.
{}From \ref{_hyperco_SLHS_Theorem_}
 it follows that the integral closure of
$\hat \calo_xM$ is also a regular ring, where 
$\hat \calo_xM$ is an adic completion of the
local ring of holomorphic functions on $(M, I)$ in a neighbourhood
of $x$. Applying \cite{_Matsumura:Commu_Alge_},
Chapter 9, Proposition 24.E again, we obtain that the
integral closure of $\calo_x M$ is regular. This proves
\ref{_normali_smooth_Corollary_}
\endproof

\hfill

\theorem \label{_desingu_Theorem_}
(Desingularization theorem)
Let $M$ be  a hypercomplex variety
$I$ an integrable induced complex structure.
Let \[ \widetilde{(M, I)}\stackrel n\arrow (M,I)\] 
be the normalization of
$(M,I)$. Then $\widetilde{(M, I)}$ is smooth and
has a natural hypercomplex structure $\c H$, such that the associated
map $n:\; \widetilde{(M, I)} \arrow (M,I)$ agrees with $\c H$.
Moreover, the hypercomplex manifold $\tilde M:= \widetilde{(M, I)}$
is independent from the choice of induced complex structure $I$.

\hfill

{\bf Proof:} The variety $\widetilde{(M, I)}$
is smooth by \ref{_normali_smooth_Corollary_}.
Let $x\in M$, and $U\subset M$ be a neighbourhood of $x$.
Let ${\goth R}_x(U)$ be the set of irreducible components of 
$U$ which contain $x$. There is a natural map
$\tau: {\goth R}_x(U) \arrow Irr(Spec \hat\calo_xM)$,
where $Irr(Spec \hat\calo_xM)$ is a set of irreducible 
components of $Spec \hat \calo_xM$, where 
$\hat \calo_xM$ is a completion of 
$\calo_xM$ in $x$. Since $\calo_x M$ is Henselian 
(\cite{_Raynaud_}, VII.4), there exist a neighbourhood $U$ of $x$
such that $\tau: {\goth R}_x(U) \arrow Irr(Spec \hat \calo_xM)$
is a bijection. Fix such an $U$. Since $M$ is a space 
with locally homogeneous singularities, the irreducible decomposition
of $U$ coinsides with the irreducible decomposition of
the tangent cone $Z_x M$. 

Let $\coprod U_i \stackrel u \arrow U$ be the morphism
mapping a disjoint union of irreducible components of $U$ 
to $U$. By \ref{_cone_flat_Theorem_}, the $x$-completion
of $\calo_{U_i}$ is regular. Shrinking $U_i$ if necessary,
we may assume that $U_i$ is smooth. Then, the morphism
$u$ coinsides with the normalization of $U$. 

For each variety $X$, we denote by $X^{ns}\subset X$ 
the set of non-singular
points of $X$. Clearly, $u(U_i) \cap U^{ns}$ is a connected
component of $U^{ns}$. Therefore, $u(U_i)$ is trianalytic
in $U$. By \ref{_triana_hyperco_Remark_}, $U_i$ has a natural
hypercomplex structure, which agrees with the map $u$.
This gives a hypercomplex structure on the normalization
$\tilde U := \coprod U_i$. Gluing these hypercomplex structures,
we obtain a hypercomplex structure $\c H$ on 
the smooth manifold$\widetilde{(M, I)}$.
Consider the normalization map $n:\; \widetilde{(M, I)} \arrow M$,
and let $\tilde M^{n}:= n^{-1}(M^{ns})$. Then, 
$n\restrict{\tilde M^{n}} \tilde M^{n}\arrow M^{ns}$ 
is an etale finite covering which is compatible with the hypercomplex
structure. Thus, $\c H\restrict{\tilde M^{n}}$ can be obtained
as a pullback from $M$. Clearly, a hypercomplex structure 
on a manifold is uniquely
defined by its restriction to an open dense subset. We
obtain that $\c H$ is independent from the choice of $I$.
\endproof

\subsection{Integrability of induced complex structures}

\theorem\label{_all_indu_comple_integra_Theorem_}
Let $M$ be a hypercomplex variety, $I$ an induced complex
structure. Then $I$ is integrable.

\hfill

{\bf Proof:}  Kaledin has proven
\ref{_all_indu_comple_integra_Theorem_} for smooth $M$
(\cite{_Kaledin_}).
Let $\tilde M$ be a desingularization of $M$, which is
hypercomplex. Then $I$ induces an integrable almost complex 
structure on $\tilde M$. From the local structure of the
singularities of $M$, is is clear that
$M$ is obtained from $\tilde M$ by gluing pairwise
certain trianalytic subvarieties $X_i\subset \tilde M$.
Since $I$ induces an integrable complex structure
on $\tilde M$, the $X_i$ are complex subvarieties
of $(\tilde M, I)$, and the identification procedures
are complex analytic with respect to $I$. Performing
these identification morphisms on $(\tilde M, I)$, we 
obtain a complex variety $M'$ such that $(M,I)$
is the underlying almost complex variety.
This proves \ref{_all_indu_comple_integra_Theorem_}.
\endproof


\section{Twistor spaces of hypercomplex varieties}
\label{_twistors_Section_}


Let $M$ be a hypercomplex variety, $M_\R$ the underlying
real analytic variety. Consider the variety
$\Tw_\R:= M_\R \times S^2$, where $S^2$ is the
2-dimensional sphere identified with the real variety
underlying $\C P^1$. We endow $\Tw_\R$ with an almost complex structure
as follows. We have a decomposition 
\[ 
   \Omega^1 \Tw_\R = \pi^* \Omega^1 S^2 \oplus \sigma^* \Omega^1 M_\R,
\]
where $\pi:\; \Tw_\R \arrow S^2$, $\sigma:\;  \Tw_\R \arrow M$ are
the natural projection maps.  Let $C$
be the natural map  
$\Omega^1 M_\R\otimes \Bbb H\stackrel C \arrow \Omega^1 M_\R$ 
arising from the quaternionic structure. 
We identify the points of $S^2$ with
induced complex structures on $M$, as usually, which
are quaternions $L\in \Bbb H$, $L^2=-1$. This gives a 
natural real analytic map $i:\; S^2 \arrow \Bbb H$.
A composition of $C$ and $i$ gives an endomorphism
$\c I_0:\; \sigma^* \Omega^1 M_\R \arrow \sigma^* \Omega^1 M_\R$.
In terms of the fibers, the endomorphism $\c I_0$ can
be described as follows. For $(s,m) \in S^2 \times M_\R = \Tw_\R$,
$\c I_0$ acts on 
$\sigma^* \Omega^1 M_\R\restrict{(s,m)}=\Omega^1 M_\R\restrict m$ by
$I_s$, where $I_s= i(s) \in \Bbb H$ is the induced complex
structure corresponding to $s\in S^2$. Since $S^2$ is identified
with $\C P^1$, this space has a natural complex structure
$I_{\C P^1}:\; \Omega^1 S^2\arrow \Omega^1 S^2$.
Let $\c I$ be an almost complex structure  
$\c I:\; \Omega^1 \Tw_\R  \arrow \Omega^1 \Tw_\R$ acting 
as $\pi^* I_{\C P^1}$ on $\pi^* \Omega^1 S^2$ and as 
$\c I_0$ on $\sigma^* \Omega^1 M_\R$. 

\hfill

\claim\label{_twi_integra_Claim_}
The constructed above almost complex structure on $\Tw_\R$ 
is integrable.

{\bf Proof:} For $M$ non-singular, this is proven by D. Kaledin
\cite{_Kaledin_}. For $M$ singular, the proof essentially
repeats the proof of \ref{_all_indu_comple_integra_Theorem_}:
we apply the desingularization 
theorem (\ref{_desingu_Theorem_}), and then Kaledin's result. \endproof

\hfill

\definition\label{_twistor_Definition_}
Let $M$ be a hypercomplex variety. Consider the 
complex variety $(\Tw, \c I)$ obtained in \ref{_twi_integra_Claim_}.
Then $\Tw$ is called {\bf a twistor space} of $M$. 

\hfill

It is possible to characterize the hypercomplex varieties
in terms of the twistor spaces. This characterization is
the main purpose of the present paper.

Consider the unique anticomplex involution 
$\iota_0:\; \C P^1 \arrow \C P^1$ with no fixed
points. This involution is obtained by central symmetry
with center in $0$ if we identify $\C P^1$ with
a unit sphere in $\R^3$. Let $\iota:\; \Tw\arrow \Tw$ 
be an involution of the twistor space mapping 
$(s,m)\in S^2 \times M = \Tw$ to $(\iota_0(s), m)$.
Clearly, $\iota$ is anticomplex.

\hfill

\definition
Let $s:\; \C P^1 \arrow \Tw$ be a section of the
natural holomorphic projection
$\pi:\; \Tw \arrow \C P^1$, $s\circ \pi = Id\restrict{\C P^1}$.
Then $s$ is called {\bf the twistor line}. The space $\Sec$
of twistor lines is finite-dimensional and equipped with
a natural complex structure, as follows from deformation theory
(\cite{_Douady_}).

\hfill

Let $\Sec^\iota$ be the space of all lines $s\in \Sec$ which are fixed
by $\iota$. The space $\Sec^\iota$ is equipped with a structure
of a real analytic space. We have a natural map 
$\tau:\; M_\R \arrow \Sec^\iota$ associating to $m\in M$ the line
$s:\; \C P^1 \arrow \Tw$, $s(x) = (x,m) \in S^2 \times M = \Tw$.
Such twistor lines are called {\bf horizontal twistor lines.}
Denote the set of horizontal twistor lines by $\Hor\subset \Sec$.

\hfill

\lemma\label{_hori-compo-inSec^iota_Claim_}
 Let $M$ be a hypercomplex variety, $\Tw$ its twistor space and
$\tau:\; M_\R \arrow \Sec^\iota$ the real analytic 
map constructed above. 
Then $\tau$ is a closed embedding identifying $M$ with one of
connected components of $\Sec^\iota$.

\hfill

{\bf Proof:} By the Desingularization Theorem (\ref{_desingu_Theorem_}),
we may assume that $M$ is smooth. For smooth $M$,
\ref{_hori-compo-inSec^iota_Claim_} is a 
well-known statement which can be easily deduced from the deformation
theory. For details, the reader is referred to \cite{_HKLR_}. \endproof

\hfill

The following data suffice to recover the hypercomplex 
variety $M$:

\begin{equation}\label{_twistor_data_Equation_}
\begin{minipage}[m]{0.8\linewidth}
\begin{itemize}
\item A complex analytic variety $\Tw$, equipped with a 
morphism $\pi:\; \Tw \arrow \C P^1$.
\item An anticomplex involution $\iota:\; \Tw \arrow \Tw$
such that $\iota\circ \pi = \pi\circ \iota_0$
\item A choice of connected component $\Hor $ of 
$\Sec^\iota\subset \Sec$.
\end{itemize}
\end{minipage}
\end{equation}

The hypercomplex variety $M$ is reconstructed as follows.
The real analytic structure on $M = \Hor $ comes from
$\Sec^\iota$. For $I \in \C P^1$, consider the map 
$p_I:\; \Hor  \arrow \pi^{-1}(I) \subset \Tw$, 
$s\in \Hor  \arrow s(I) \in \Tw$. This identifies 
$\Hor $ with $\pi^{-1}(I) \subset \Tw$.
We obtained a complex structure $I$ on $M$,
for each $I\in \C P^1$. Identifying $\C P^1$ with
a subset of quaternions, we recover the original quaternion action
on $\Omega^1 M$. The data \eqref{_twistor_data_Equation_}
satisfies the following properties (condition (ii) is
implicit in the quaternionic action). 

\begin{equation}\label{_twistor_properties_Equation_}
\begin{minipage}[m]{0.8\linewidth}
\begin{description}
\item[(i)] For each point $x\in \Tw$, there is a unique line 
$s\in \Hor \subset \Sec^\iota$, passing through $x$. 
\\[2mm]
\begin{minipage}[m]{0.8\linewidth}
{\it For every $s\in \Sec$, we identify the
image of $s$, $\im s \subset \Tw$, with $\C P^1$.}
\end{minipage}
\item[(ii)] For every line $s\in \Hor\subset \Sec^\iota$,
the conormal sheaf
\[ N^*_s = \ker\left( \Omega^1 \Tw\restrict {\im s} \stackrel {s^*}\arrow
   \Omega^1 (\im s)\right)
\]
of $\im s$ is isomorphic to 
$\calo(-1) \oplus \calo(-1) \oplus ... \oplus \calo(-1)$.
\end{description}
\end{minipage}
\end{equation}

\hfill

\definition \label{_twi_hyperco_type_Definition_}
The data \eqref{_twistor_data_Equation_}
satisfying the conditions (i), (ii) of 
\eqref{_twistor_properties_Equation_}
are called {\bf a twistor space of hypercomplex type}.
We have shown how to associate a twistor space of 
hypercomplex type to every hypercomplex variety. 

Denote the corresponding functor by $\c F$.

\hfill

Condition (ii) of \eqref{_twistor_properties_Equation_}
can be replaced by the following condition.

\begin{equation}\tag{\ref{_twistor_properties_Equation_}}
\begin{minipage}[m]{0.8\linewidth}
\begin{description}
\item[(ii$'$)] For every line $s\in \Hor\subset \Sec^\iota$,
there exists an open neighbourhood $U\subset \Tw$ of $\im s$,
such that for every $x, y \in U$, $\pi(x) \neq \pi(y)$,
there exists a unique twistor line $s_{x,y}$ passing through
$x$ and $y$, provided that $x$ and $y$ belong to the same 
irreducible component of $U$.
\end{description}
\end{minipage}
\end{equation}

Condition (ii) of \eqref{_twistor_properties_Equation_}
should be thought of as a linearization of 
\eqref{_twistor_properties_Equation_} (ii$'$).
In the subsequent section, we shall see that these
conditions are equivalent.

\hfill

\definition\label{_Deligne_Simpson_Definition_}
The data \eqref{_twistor_data_Equation_}
satisfying the conditions (i), (ii$'$) of 
\eqref{_twistor_properties_Equation_}
are called {\bf a twistor space of Deligne-Simpson type}.
These conditions were proposed by Deligne and Simpson
(\cite{_Simpson:hyperka-defi_}, \cite{_Deligne:defi_})
in order to define singular hyperk\"ahler manifolds.


\section{Twistor spaces of Deligne-Simpson type}
\label{_Deli_Si_equi_hyperco_Section_}


The main result of this section is the following theorem.

\hfill

\theorem \label{_Deli_Simpsi_equi_infinite_Theorem_}
Let $(\Tw, \pi, \iota, \Hor)$ be the data
of \eqref{_twistor_data_Equation_} satisfying
condition (i) of \eqref{_twistor_properties_Equation_}.
Then the conditions (ii) and (ii$'$) are equivalent.
In other words, $(\Tw, \pi, \iota, \Hor)$ is a twistor
space of hypercomplex type if and only if $(\Tw, \pi, \iota, \Hor)$
is a twistor space of Deligne--Simpson type.

\hfill

The proof of \ref{_Deli_Simpsi_equi_infinite_Theorem_}
takes the rest of this section.

\hfill

Under assumptions of \eqref{_twistor_data_Equation_},
\eqref{_twistor_properties_Equation_} (i)
consider the map $\sigma:\; \Tw \arrow \Hor$ associating
to a point $x\in \Tw$ the unique horisontal line
passing through this point. This map is continuous,
and induces a homeomorphism 
$\sigma \times \pi:\; \Tw \arrow \Hor \times \C P^1$.

\hfill

\lemma\label{_Deli_Si_local_Lemma_} 
Let $(\Tw, \pi, \iota, \Hor)$ be the data
of \eqref{_twistor_data_Equation_} satisfying
condition (i) of \eqref{_twistor_properties_Equation_},
and $U\subset \Hor$ an arbitrary open subset. Then 
$\sigma^{-1}(U)$ is preserved by $\iota$.

{\bf Proof:} Let $s \in U \subset \Hor$. Then $\iota$ preserves a line
$\im s \subset \Tw$, and thus, 
$\iota(\im s) \subset \sigma^{-1}(s) \subset \sigma^{-1}(U)$.
This proves \ref{_Deli_Si_local_Lemma_}. \endproof

\hfill

We prove the implication ``$(\Tw, \pi, \iota, \Hor)$
of hypercomplex type'' $\Rightarrow$  ``$(\Tw, \pi, \iota, \Hor)$
of Deligne--Simpson type''.
The statement of \ref{_Deli_Simpsi_equi_infinite_Theorem_}
is local in $\Hor$, as follows from
\ref{_Deli_Si_local_Lemma_}. Consider the evaluation maps 
$p_I:\; \Hor \arrow \pi^{-1}(I)$, $s \arrow s(I)$, defined
for all $I\in \C P^1$. The map $\sigma$ gives a homeomorphism
$\pi^{-1}(I) \stackrel \sigma \arrow \Hor$. A homeomorphism
preserves the dimension of the variety. Thus, for a smooth
point $x\in \Tw$, the point $\sigma(x)\in \Hor$ is
equidimensional
\footnote{Equidimensional point of $X$ is a point where
all irreducible components of $X$ have the same dimension.}
 in $\Hor$. Using the homeomorphism
$\sigma \times \pi:\; \Tw \arrow \Hor \times \C P^1$,
we find that for a smooth point $x\in \Tw$ 
and $y\in \Tw$ satisfying $\sigma(y)=\sigma(x)$, 
the point $y$ is equidimensional in $\Tw$.

The real dimensions of the varieties
$\Hor$ and $\pi^{-1}(I)$ are equal. On the other hand,
for every $m\in \Hor$, the dimension of the tangent space
$T_{p_I(m)}\pi^{-1}(I)$ is equal to the dimension of
$N^*_m\restrict{p_I(m)}$. Since $N^*_m\restrict{p_I(m)}$
is a bundle, the dimension of $T_{p_I(m)}\pi^{-1}(I)$
is the same for all $I\in \C P^1$.
The local ring is regular if and only if the dimension
of its tangent space is equal to the dimension of the ring.
Dimensions of the varieties
$\pi^{-1}(J)$ and $\pi^{-1}(I)$, and the dimensions of the
corresponding tangent spaces are equal.
Thus, for a smooth point $m\in \pi^{-1}(I)$, and every $J\in \C P^1$,
the points $p_J(\sigma (m))\in \pi^{-1}(J)$ are
smooth in $\pi^{-1}(J)$. 

We obtained the following result.

\hfill

\claim \label{_hori_smooth_Claim_}
Let $(\Tw, \pi, \iota, \Hor)$ be a twistor space 
of hypercomplex type, and $m\in \Tw$ a smooth point.
Let $\sigma$ denote the natural continuous map 
$\sigma:\; \Tw \arrow \Hor$. Then, for every point $m'\in \Tw$ such that
$\sigma(m) = \sigma(m')$, $m'$ is smooth. 

\endproof

\hfill

\lemma\label{_hype+smoo=>Deli-Simps_Lemma_}
Let $(\Tw, \pi, \iota, \Hor)$ be a twistor space 
of hypercomplex type, and $m\in \Tw$ a smooth point.
Then 
\begin{description}
\item[(i)]
$\sigma(m)$ is a smooth point of $\Hor$.
\item[(ii)]
Moreover, for a smooth neighbourhood $U$ of $m\in Tw$, 
$(\sigma^{-1}(\sigma(U))\subset \Tw, \pi, \iota, \sigma(U))$ 
is a twistor space of Deligne-Simpson type.
\end{description}

{\bf Proof:} Let $s= \im m\subset \Tw$ be the horisontal twistor line
corresponding to $m$. From the deformation theory we know that
the deformations of a smooth curve $s$ are classified by the sections of
the normal bundle $\Gamma(Ns)$, with obstructions corresponding
to $H^1(Ns)$. The cohomology space $H^1(Ns)$ vanishes, because
$Ns=\bigoplus \calo(1)$ is ample. Thus, $\Gamma(Ns)$ has a dimension
$2 (\dim \Tw-1)$. For a small deformations $s'$ of $s$,
$\Gamma(Ns')= \bigoplus \calo(1)$, since $\bigoplus \calo(1)$
is semistable. Thus, $T_{s'}\Sec$ is constant in a neighbourhood of $s$,
where $\Sec$ is the space of sections of $\pi:\; \Tw\arrow \C P^1$.
This implies that $s$ is a smooth point of $\Sec$, and hence, 
$s$ is a smooth point of $\Hor$. \ref{_hype+smoo=>Deli-Simps_Lemma_} 
(i) is proven.

To prove (ii), let $I$, $J$ be distinct points in $\C P^1$
and consider the map 
$p_{IJ}:\; \Sec \arrow \pi^{-1}(I) \times \pi^{-1}(J)$,
$\gamma \arrow (\gamma(I), \gamma(J))$.
We have to prove that $p_{IJ}$ has invertible
differential in $s$. The tangent space $T_s \Sec$ is, as we have seen,
 $\Gamma(Ns)$. The differential of the map
$p_I:\; \Sec\arrow \pi^{-1}(I)$, $\gamma \arrow \gamma(I)$
coinsides with the restriction map 
$r_I:\;\Gamma(Ns) \arrow Ns\restrict{I}$.
Since $Ns$ is $\bigoplus \calo(1)$, the differential
$dp_{IJ} = r_I\times r_J$ is an isomorphism (a section of
$\bigoplus \calo(1)$ is uniquely determined by
its value in two distinct points). This proves 
\ref{_hype+smoo=>Deli-Simps_Lemma_} (ii). \endproof

\hfill

We return to the proof of an implication ``$(\Tw, \pi, \iota, \Hor)$
of hypercomplex type'' $\Rightarrow$  ``$(\Tw, \pi, \iota, \Hor)$
of Deligne--Simpson type''. 
Let $\Tw^{ns}$ be the set of non-singular
points of $\Tw$, $I$, $J$ be two distinct points of
$\C P^1$. Let $W\subset \Sec$ be an open neighbourhood 
of a horisontal line $s\in \Sec$, such that its closure
$\bar W$ is compact and $\bar V_{IJ}$ the set of all triples 
$(x,y, s_{xy}) \in \pi^{-1}(I)\times \pi^{-1}(J)\times \Sec$
such that  $s_{xy}\in \bar W$, and $s_{xy}$ is a twistor line
passing through $x$ and $y$. Let 
$V_{IJ}\subset \bar V_{IJ}$  be the set of the triples 
$(x,y, s_{xy})\in \bar V_{IJ}$, for which 
the corresponding twistor line $s_{xy}$ belongs
to the non-singular part of $\Tw$.

\hfill

\lemma \label{_closure_of_lines_through_ns_Lemma_}
In the above assumptions, consider the 
forgetful map 
$p:\; \bar V_{IJ} \arrow \pi^{-1}(I)\times \pi^{-1}(J)$.
Then  $p\left(\bar V_{IJ}\right)$
is a closure of $p\left( V_{IJ}\right)$.

\hfill

{\bf Proof:} The space $\bar V_{IJ}$ is compact, and hence
its image is compact and thus closed. It remains to show that
$V_{IJ}$ is dense in $\bar V_{IJ}$. By
\ref{_hori_smooth_Claim_}, for each smooth point
$m\in \Tw$, a neighbourhoor of $\sigma(m)\in \Sec$
belings to $V_{IJ}$. Since the set of smooth points
of $\Tw$ is dense in $\Tw$, for $\bar V_{IJ}$ sufficiently
small, $V_{IJ}$ is dense in $\bar V_{IJ}$.
\endproof

\hfill

Let $U$ be an open subset of $\Tw$, and
$X_U\subset \pi^{-1}(I)\times \pi^{-1}(J)$ be the set of
all $(x,y) \in \pi^{-1}(I)\times \pi^{-1}(J)$ belonging to the
same irreducible component of $U\subset \Tw$. 
Consider the forgetful map 
$p:\; \bar V_{IJ} \arrow  \pi^{-1}(I)\times \pi^{-1}(J)$.
Clearly, the image of $p$ intersected with $U\times U$
lies in $X_U$, so we may assume that
$p$ maps $\bar V_{IJ}$ to $X_U$.
Computing the differential of $p$ as in the proof of
\ref{_hype+smoo=>Deli-Simps_Lemma_}, we find that
$dp$ is locally injective for $s_{xy}$ in a neighbourhood
of $\Hor$. To prove the condition of
Deligne and Simpson, it remains to show that
$p$ is locally a surjection onto $X_U$, for sufficiently
small $U$. By \ref{_hype+smoo=>Deli-Simps_Lemma_},
the image of $p\restrict{V_{IJ}}$ 
is dense in $X_U\cap \Tw^{ns}\times \Tw^{ns}$,
for $U$ sufficiently small. 
On the other hand, the closure of  
$\im p\restrict{V_{IJ}}$ 
is $\im p\restrict{\bar V_{IJ}}$
by \ref{_closure_of_lines_through_ns_Lemma_},
so $p$ is locally a surjection.
We proved that the Deligne-Simpson condition 
holds for all twistor spaces of hypercomplex type.

Assume now that $(\Tw, \pi, \iota, \Hor)$ is a twistor space 
of Deligne--Simpson type. Let $s\in \Hor$. 
Consider the 
evaluation map $p_{IJ}:\; \Sec \arrow \pi^{-1}(I)\times \pi^{-1}(J)$,
$I\neq J\in \C P^1$. By Deligne-Simpson's condition,
$p_{IJ}$ induces an isomorphism 
\begin{equation}\label{_diffe_eva_Tsec_Equation_}
   dp_{IJ}:\; T_s \Sec 
   \arrow T_{s(I)}\pi^{-1}(I)\times T_{s(I)}\pi^{-1}(I)
\end{equation}
of the tangent spaces. Thus, dimension of $T_{s(I)} \pi^{-1}(I)$
is independent from the choice of $I$. We obtain that
the conormal sheaf $N^*s$ is a bundle, and it makes sense
to speak of the normal bundle $Ns$.

Through each point in a neighbourhood of $s$ passes
a deformation of $s$. Thus, $T_s \Sec= \Gamma(Ns)$; there is no
first order obstructions to the deformation.
 Since $T_s \Sec= Ns$, the map 
\eqref{_diffe_eva_Tsec_Equation_} can be interpreted as
\[ 
 dp_{IJ}:\;\Gamma(Ns) \arrow Ns\restrict{I} \times Ns\restrict{J}
\]
Since \eqref{_diffe_eva_Tsec_Equation_} is an isomorphism,
$Ns$ is a bundle which is isomorphic to $\bigoplus \calo(1)$.
\ref{_Deli_Simpsi_equi_infinite_Theorem_} is proven. 
\endproof


\section{Hypercomplex varieties and twistor spaces of hypercomplex
type.}
\label{_hype_type_equi_hype_Section_}


The main result of this paper is the following theorem.

\hfill

\theorem \label{_hype_vari_and_twi_equiva_Theorem_}
Consider the functor $\c F$ of \ref{_twi_hyperco_type_Definition_},
associating to a hypercomplex variety the corresponding 
twistor space of hypercomplex type.
Then $\c F$ is equivalence of categories.

\hfill

{\bf Proof:} We have shown how to recover the hypercomplex
structure from the twistor space. This proves  that $\c F$ 
is full and faithful. It remains to show that each twistor space
of hypercomplex type is obtained from a hypercomplex variety. 
Thus, \ref{_hype_vari_and_twi_equiva_Theorem_} is implied
by the following statement.

\hfill

{\bf \ref{_hype_vari_and_twi_equiva_Theorem_}$'$}
Let \[ (\Tw, \pi, \iota, \Hor)\] be a twistor space 
of hypercomplex type. Then $\Hor$ admits a hypercomplex 
structure $\c H$,
such that $(\Tw, \pi, \iota, \Hor)$ is a twistor space 
of $(\Hor, \c H)$. 

\hfill

The rest of this section is taken by the proof of
\ref{_hype_vari_and_twi_equiva_Theorem_}$'$. 

\hfill

\lemma\label{_irre_compo_twi_hype_type_Lemma_}
Let $\Hor = \bigcup H_i$ be an irreducible decomposition
of the variety $\Hor$. Let $\Tw_i := \sigma^{-1}(H_i)$,
where $\sigma:\; \Tw\arrow \Hor$ is the standard 
continuous map. Then $\Tw_i$ is preserved by $\iota$, and
\[ \left( Tw_i, \pi\restrict{\Tw_i}, 
   \iota\restrict{\Tw_i}, H_i\right)
\]
is a twistor space of hypercomplex type. 

\hfill

{\bf Proof:}
It is clear from the definition that $\Tw_i$ is preserved by
$\iota$ and that for every $m\in \Tw_i$ there is a unique
horisontal line $s\in H_i$ passing through $m$. It remains to prove the
condition (ii) of \eqref{_twistor_properties_Equation_},
or, equivalently, condition (ii)$'$ 
of \eqref{_twistor_properties_Equation_}.
But, since $\Tw$ satisfies the condition (ii)$'$ 
of \eqref{_twistor_properties_Equation_}, $\Tw_i$ 
satisfies this condition automatically. \endproof

\hfill

\claim \label{_homeomo_maos_irre_compo_Claim_}
Let $f:\; X \arrow Y$ be a homeomorphism of complex analytic
varieties. Then $f$ maps irreducible components of $X$ to irreducible 
components of $Y$.

{\bf Proof:} Clear. \endproof 

\hfill

\lemma \label{_p_I_iso_Lemma_}
Let $(\Tw, \pi, \iota, \Hor)$ be a twistor space 
of hypercomplex type, and $I\in \C P^1$. Consider the map
$p_I:\; \Hor \arrow \pi^{-1}(I)$, $s \arrow s(I)$.
Let $p:\; \Hor \times \C P^1\arrow \Tw$ map $(s, I)$ to $s(I)$.
Then $p$ and $p_I$ induce isomorphisms of corresponding
real analytic varieties.

\hfill

{\bf Proof:} Using \ref{_homeomo_maos_irre_compo_Claim_},
\ref{_irre_compo_twi_hype_type_Lemma_}
and \ref{_Deli_Si_local_Lemma_}, we may assume that the fiber
$\pi^{-1}(I)$ is locally irreducible and thus, equidimensional. 
The real dimensions of the varieties $\Hor$, $\pi^{-1}(I)$
are clearly equal. Thus, to show that $p$, $p_I$ induce
isomorphisms, it suffices to show that corresponging maps
of local rings are surjective. By Nakayama, for this
we need to show that $p$, $p_I$ induce surjection on the 
Zariski tangent spaces. The differential of the evaluation
map $ev_{I}:\; \Sec \arrow \pi^{-1}(I)$, $s\arrow s(I)$
is the standard restriction map
$r_I:\; \gamma \arrow \gamma\restrict{I}$, where 
$\gamma\in T_s \Sec = \Gamma(Ns)$, and 
$ \gamma\restrict{I}\in Ns\restrict{I} = T_{s(I)}\pi^{-1}(I)$.
Thus, $dp_I$ is a composition of $r_I$ and the embedding
$T_s \Hor \hookrightarrow T_s\Sec$. The image
of $T_s \Hor \hookrightarrow T_s\Sec$ coinsides with the
set of fixed points of the involution 
$d\iota:\; T_s\Sec \arrow T_s \Sec$. Thus,
$dp_I$ is an isomorphism by the following
trivial result of linear algebra.

\hfill

\lemma \label{_tota_real_iso_linear_alge_Lemma_}
Let $V$ be a vector space of complex dimension $2n$,
$W$  a vector space of complex dimension $n$, and
$\phi:\; V \arrow W$ an epimorphism. Let $V'\subset V$
be a totally real subspace of real dimension $2n$.
Then $\phi\restrict{V'}:\; V' \arrow W$ is an isomorphism.

\endproof

A similar argument proves that $p$ is also an isomorphism.
\ref{_p_I_iso_Lemma_} is proven. \endproof

\hfill

We obtained that the real analytic variety $\Hor$ is isomorphic 
to one underlying $\pi^{-1}(I)$, for all $I\in \C P^1$. This
gives a set of integrable almost complex structures on
$\Hor$, parametrized by $\C P^1$. The following 
linear algebraic argument shows that these complex structures
satisfy quaternionic relations. This finishes the proof of 
\ref{_hype_vari_and_twi_equiva_Theorem_}.

\hfill

Let $F$ be an $n$-dimensional holomorphic vector bundle on
$\C P^1$, $F \cong \oplus \calo(1)$. Consider
the restriction maps $r_I:\; \Gamma(F) \arrow F\restrict{I}$, defined
for each $I\in \C P^1$. Let $W\subset \Gamma(F)$ be a totally
real subspace of real dimension $2n$. By
\ref{_tota_real_iso_linear_alge_Lemma_}, 
$r_I$ induces an isomorphism between
$F\restrict{I}$ and $W$. This gives a
set of complex structures on $W$, parametrized by 
$\C P^1$.

\hfill

\claim
These complex structures satisfy quaternionic relations.

{\bf Proof:} Clear. \endproof


\section{Some applications}
\label{_twi_applications_Section_}


\subsection{Hypercomplex spaces}

Using \ref{_hype_vari_and_twi_equiva_Theorem_}, it is
possible to generalize the definition of hypercomplex varieties,
allowing nilpotents, in such a way that a reduction of a 
hypercomplex space is a hypercomplex variety.

Consider the anticomplex 
involution $\iota_0:\; \C P^1 \arrow \C P^1$ defined in 
\ref{_twistors_Section_}. 

\hfill

\definition\label{_hyperco_spaces_Definition_}
(Hypercomplex spaces)
Let $\Tw$ be a complex analytic space,
$\pi:\; \Tw\arrow \C P^1$ a holomorphic map,
and $\iota:\; \Tw \arrow \Tw$ an anticomplex automorphism,
such that $\iota\circ\pi = \pi\circ \iota_0$. Let $\Sec$ be the
space of sections of $\pi$ equipped with a structure of a complex
analytic space, and $\Sec^\iota$ be the real analytic space of 
sections $s$ of $\pi$ satisfying $s\circ \iota_0 = \iota\circ s$. Let $\Hor$
be a connected component of $\Sec^\iota$. Then $(\Tw, \pi,\iota,\Hor)$
is called {\bf a hypercomplex space} if
\begin{description}
\item[(i)] For each point $x\in \Tw^r$, there exist a unique line
$s\in \Hor^r$ passing through $x$, where $\Tw^r$, $\Hor^r$ is
a reduction of $\Tw$, $\Hor$.

\item[(ii)] Let $s\in \Hor$, and $U\subset \Tw$ be a neighbourhood of
$s$ such that an irreducible decomposition of $U$ coinsides
with the irreducible decomposition of $\Tw$ in a neighbourhood
of $s\subset \Tw^r$. Let 
\[ \bar X:= \pi^{-1}(I)\times \pi^{-1}(J)\cap
   U\times U,
\] where $I$, $J$ distinct points of $\C P^1$.  Let 
$p_{IJ}:\; U \arrow \bar X\subset \pi^{-1}(I)\times \pi^{-1}(J)$
be the evaluation map, $s\arrow (s(I), s(J))$. Then there
exist a closed subspace $X\subset \bar X$, obtained as a union
of some of irreductible components of $\bar X$, and an open
neighbourhood $V\subset \Sec$ of $s\in \Sec$, such that
$p_{IJ}$ is an open embedding of $V$ to $X$.
\end{description}


\subsection{Stable bundles over hyperk\"ahler manifolds}
\label{_hyperholomo_Subsection_}


Let $M$ be a compact hyperk\"ahler manifold, 
$I$ an induced complex structure and $B$ a stable holomorphic
bundle over $(M, I)$, such that the first two Chern classes 
$c_1(B)$, $c_2(B)$ are $SU(2)$-invariant, with respect to the
natural action of the group $SU(2)$ on the cohomology of $M$ 
(\ref{_SU(2)_commu_Laplace_Lemma_}).

Recall that $SU(2)$ acts
on the space of differential forms on $M$.
This allows us to speak of $SU(2)$-invariant differential forms,
for instance of connections with $SU(2)$-invariant curvature.
In \cite{_Verbitsky:Hyperholo_bundles_}, the following theorem
was proven.

\hfill

\theorem \label{_hyperholomo_conne_Theorem_}
There exist a unique Hermitian connection $\nabla$ on $B$
such that its curvature $\Theta$ is $SU(2)$-invariant.
Conversely, if such connection exitst on a holomorphic
bundle $B$ over $(M,I)$, then $B$ is a direct 
sum of stable bundles with $SU(2)$-invariant 
Chern classes.

\endproof

\hfill

We show that the moduli space $\Def(B)$ of deformations of 
$B$ is hypercomplex. Consider the twistor space $\Tw$ of $M$,
and a standard real analytic map $\sigma:\; \Tw \arrow M$.
Let $\sigma^* B$ be the pullback of $B$ equipped with the
connection which is trivial along the fibers of $\sigma$. 
{}From twistor transform (\cite{_NHYM_}, Section 5) is clear that
$\sigma^*B$ is holomorphic. Restricting $\sigma^*B$ to the
fibers of $\pi:\; \Tw \arrow \C P^1$, we obtain holomorphic
bundles $B_J$ on $\pi^{-1}(J) = (M, J)$. By 
\ref{_hyperholomo_conne_Theorem_}, $B_J$ is stable. Let $\hat{\Tw}$
be the moduli space of sheaves of type $i^J_* F$, where
$i^J:\; (M, J)=\pi^{-1}(J)\hookrightarrow \Tw$ is the
natural embedding, and $F$ a stable holomorphic
bundle which is a deformation of $B_J$ (i. e. belongs in the
same deformation class). Then $\hat{\Tw}$ is equipped with a
holomorphic fibration $\hat \pi:\; \hat{\Tw}\arrow \C P^1$, 
$i^J_* F\arrow J$. Mapping $F$ to $\iota^* F$, we obtain 
an anticomplex involution $\hat \iota$ of $\hat{\Tw}$. The operation
$(B,J)\arrow B_J$ gives an $\hat\iota$-invariant section of $\hat\pi$,
parametrized by $\Def(B)$. To show that thus obtained
quadruple $(\hat{\Tw}, \hat \pi, \hat\iota, \Def(B))$ is
hypercomplex, it suffices to prove the condition (ii) of
\eqref{_twistor_properties_Equation_}.
Equivalently, we may prove (ii$'$) of 
\ref{_Deligne_Simpson_Definition_}. On the other hand, 
Proposition 2.19 of \cite{_NHYM_} implies (ii$'$).
This gives another proof that the space of stable deformations
of $B$ is hypercomplex, in addition to that given in
\cite{_Verbitsky:Hyperholo_bundles_}. 

\subsection{Quotients of hypercomplex varieties by an action
of a finite group}

Let $M$ be a hypercomplex variety and $G$ a finite group acting on $M$,
generically free. Assume that $G$ preserves the hypercomplex structure,
and acts freely outside of nonempty
finite set of fixed points, denoted by $\Gamma$.
Clearly, $\bigg(M\backslash \{\Gamma\}\bigg)/G$ 
has a natural structure of a 
hypercomplex variety.

\hfill

\proposition \label{_quotie_not_hyperco_Proposition_}
The hypercomplex structure on 
$\bigg(M\backslash \{\Gamma\}\bigg)/G\subset M/G$
cannot be extended to $M/G$.

\hfill

{\bf Proof:} 
Consider the space $\Tw/G$ fibered over $\C P^1$, with 
corresponding action of $\iota$. Then $\Hor/G$ gives an
open subset in the space of $\iota$-invariant sections
of $\pi:\; \Tw/G \arrow \C P^1$. Let $p:\; M \arrow M/G$,
$p:\; \Tw \arrow \Tw/G$ be the natural quotient maps.
If the
hypercomplex structure on $M\backslash \{\Gamma\}/G\subset M/G$
were extended to $M/G$, the space of horisontal sections 
would have been $p(\Hor)$. Applying
\ref{_hype_vari_and_twi_equiva_Theorem_}, we obtain the
following assertion. 

\hfill

\claim
The hypercomplex structure on
$M\backslash \{\Gamma\}$ is extended to $M/G$ 
if an only if for all $s:\; \C P^1 \arrow \Tw/G$, 
$s \in p(\Hor)$, the conormal sheaf of $\im s$ in
$\Tw/G$ is $\oplus \calo(-1)$. 

\endproof

\hfill

Let $s\in \Hor$
be a horisontal twistor line in $\Tw$ which 
passes through fixed point of $G$-action.
Consider its formal neighbourhood in $\Tw$.
Let $\calo(s)$ be the corresponding complete
ring over $\C P^1$ and $\calo(s)_{gr}$ be the associated
graded ring. Then the ring $\calo(s)_{gr}$ is isomorphic
to $\oplus S^*(\calo(-1))/D$, where $D$ is a graded ideal
lying in 
\[ \bigoplus\limits_{k=2}^{\infty} S^k(\calo(-1)) \]
Let $\calo(\hat s)_{gr}\subset\calo(s)_{gr} $ be the sheaf of $G$-invariant
sections of $\oplus S^*(\calo(-1))/D$. Clearly,
$\calo(\hat s)_{gr}$ is the graded ring of the 
formal neighbourhood of $\hat s \subset \Tw/G$. 
Since $x$ is an isolated fixed point of $G$ action, the group acts
on the Zariski tangent space $T_xM$ without invariants.
Therefore, $\calo(\hat s)_{gr}$ lies in 
\[ \bigoplus\limits_{k=2}^{\infty} S^k(\calo(-1))/D 
   \subset \oplus S^*(\calo(-1))/D = \calo(s)_{gr}.
\]
Thus, the conormal sheaf of $\calo(\hat s)_{gr}$
is isomorphic to $\calo(i_1) \oplus \calo(i_2) \oplus ...$,
where $i_1, ... , i_k<-2$. Therefore, $M/G$ cannot be hypercomplex.
\ref{_quotie_not_hyperco_Proposition_} is proven. \endproof

\hfill

{\bf Acknowledegments:} D. Kaledin and T. Pantev explained me
the Deligne and Simpson's 
definition of a hyperk\"ahler variety. 
They were also very helpful in correcting the errors of the
manuscript. P. Deligne kindly 
pointed out persistent errors in the presentation
of \cite{_Verbitsky:Desingu_}, trying to explain to me the
theory of real analytic spaces.
I am grateful to these and also to A. Beilinson, R. Bezrukavnikov,
M. Entov, D. Kazhdan, M. Kontsevich, 
A. Todorov and S.-T. Yau for valuable discussions.


\begin{thebibliography}{666}
 

\bibitem[A]{_Arapura_} Arapura, D.,
{\em Geometry of cohomology support loci II: integrability of
Hitchin's map}, alg-geom/9701014, 17 pages, AMS LaTeX.

\bibitem[Bea]{_Beauville_} 
Beauville, A. {\em 
Varietes K\"ahleriennes dont la premi\`ere classe de Chern est
nulle.} // J. Diff. Geom. 18, p. 755-782 (1983).


\bibitem[Bes]{_Besse:Einst_Manifo_} 
Besse, A., {\em Einstein Manifolds}, Springer-Verlag, New York (1987)

\bibitem[C]{_Calabi_} Calabi,  E.,
{\em Metriques k\"ahleriennes et fibr\`es holomorphes}, 
Ann. Ecol. Norm. Sup. {\bf 12} (1979), 269-294.  

\bibitem[De]{_Deligne:defi_} 
Deligne, P.,
{\it A letter to Carlos Simpson}, 
referred to in \cite{_Simpson:hyperka-defi_}.

\bibitem[Do]{_Douady_} 
Douady, A., {\em Le probleme des modules pour
les sous-espaces analitiques compactes d'un espace analitique
donn\'e}, Ann. Inst. Fourier, {\bf 16} (1966), vol. 1, pp. 1-95.

\bibitem[GMT]{_GMT_} Guaraldo, F., Macri, P., Tancredi, A.,
{\em Topics on real analytic spaces},
 Advanced lectures in mathematics, Braunschweig: F. Vieweg, 1986.



\bibitem[H]{_Hartshorne:Alg_Geom_} 
 Hartshorne, R., {\em Algebraic geometry}, 
 Graduate texts in mathematics, vol. 52,
 New York : Springer-Verlag, 1977.
 

\bibitem[HKLR]{_HKLR_} 
N.J. Hitchin, A. Karlhede, U. Lindstr\"{o}m,
M. Ro\v{c}ek, {\em Hyperk\"ahler metrics and supersymmetry},
Comm. Math. Phys (1987). 


\bibitem[K]{_Kaledin_} Kaledin, D.,
{\em Integrability of the twistor space for a hypercomplex manifold,}
alg-geom/9612016, 9 pages, Latex2e.

\bibitem[KV]{_NHYM_} 
Kaledin, D., Verbitsky, M.,
{\it Non-Hermitian Yang-Mills connections}, 
alg-geom 9606019 (1996), 48 pages, LaTeX 2e.


\bibitem[M] {_Matsumura:Commu_Alge_}
 Matsumura, H., {\em Commutative algebra,}  
 Mathematics lecture note series, vol. 56,
 Benjamin/Cummings Pub. Co., Reading, Mass.,  1980.
 

\bibitem[R]{_Raynaud_} Raynaud, M. {\em Anneaux Locaux Hens\'eliens}, 
Springer LNM 169, 1970.



\bibitem[S]{_Simpson:hyperka-defi_} 
Simpson, C. T.,
{\it Nonabelian Hodge theory.}
Proceedings of the International Congress of Mathematicians,
(Kyoto, 1990), 747--756, 
Math. Soc. Japan, Tokyo, 1991. 

\bibitem[V1]{_Verbitsky:Hyperholo_bundles_} 
Verbitsky M., 
{\em Hyperholomorphic bundles over a hyperk\"ahler manifold}, 
alg-geom electronic preprint 9307008 (1993), 43 pages, LaTeX,\\
also published in: 
Journ. of Alg. Geom., {\bf 5} no. 4 (1996) pp. 633-669.


\bibitem[V2]{_Verbitsky:Symplectic_II_} 
Verbitsky M., {\em Hyperk\"ahler embeddings and holomorphic 
symplectic geometry II,} alg-geom electronic preprint 9403006 (1994),
14 pages, LaTeX,
also published in: GAFA {\bf 5} no. 1 (1995), 92-104.

\bibitem[V3]{_Verbitsky:Deforma_}
Verbitsky M., {\em Deformations of trianalytic subvarieties of
hy\-per\-k\"ah\-ler manifolds}, alg-geom electronic preprint 9610010
(1996), 51 pages, LaTeX2e.

\bibitem[V4] {_Verbitsky:New_} 
Verbitsky, M. {\it New examples
of compact hyperk\"ahler manifolds}, in preparation.

\bibitem[V-d]{_Verbitsky:Desingu_}
Verbitsky M.,
{\em Desingularization of singular hyperk\"ahler varieties I,}
electronic preprint 9611015 (1996), 13 pages, LaTeX 2e
(to appear in Math. Res. Let.)

\bibitem[V-d2]{_Verbitsky:DesinguII_}
Verbitsky M.,
{\em Desingularization of singular hyperk\"ahler varieties II,}
electronic preprint 9612013 (1996), 15 pages, LaTeX 2e.


\end{thebibliography}
\end{document}